\documentclass{amsart}
\usepackage{amsmath}
\usepackage{amscd}
\usepackage{amsfonts}
\usepackage{amssymb}
\usepackage{latexsym}
\usepackage{color}
\newtheorem{theorem}{Theorem} 
\newtheorem{lemma}{Lemma}[section]     
\newtheorem{corollary}{Corollary}
\newtheorem{proposition}{Proposition}
\newtheorem{definition}{Definition}[section]
\newtheorem{remark}{Remark}
\newtheorem{example}{Example}[section]

\newtheorem{problem}{Problem}
\newtheorem{condition}{Condition}[section]
\author{T.~Nishimura}
\address{
Research Group of Mathematical Sciences, Research Institute of
Environment and Information Sciences, Yokohama National University,
Yokohama 240-8501, JAPAN. } \email{nishimura-takashi-yx@ynu.jp}
\thanks{T.~Nishimura is partially supported
by JSPS and CAPES under the Japan--Brazil research cooperative
program.}
\author{R.~Oset~Sinha}
\address
{Department of Mathematics,   Federal University of S\~ao Carlos,
S\~ao Carlos, SP, Brazil} \email{Raul.Oset@uv.es}
\thanks{\textcolor{black} {R. Oset Sinha is partially supported by FAPESP grant
no. 2013/02381-1 and DGCYT and FEDER grant no. MTM2012-33073.}}
\author{M.~A.~S.~Ruas}
\address{ICMC,
University of S\~ao Paulo, S\~ao Carlos, SP, Brazil}
\email{maasruas@icmc.usp.br}
\thanks{\textcolor{black} {M. A. S. Ruas and R. Wik Atique are partially supported by FAPESP grant no.
2014/00304-2. M. A. S. Ruas is partially supported by CNPq grant no.
303774/2008-8.}}
\author{R.~Wik Atique}
\address{ICMC,
University of S\~ao Paulo, S\~ao Carlos, SP, Brazil}
\email{rwik@icmc.usp.br}
\title 
[Liftable vector fields over corank one multigerms]
{Liftable vector fields  \\
over corank one multigerms
}
\begin{document}
\begin{abstract}
In this paper, a systematic method is given to construct all
liftable vector fields over \textcolor{black}{an analytic} multigerm
$f: (\mathbb{K}^n, S)\to (\mathbb{K}^p,0)$ \textcolor{black}
of corank at most one admitting a one-parameter stable unfolding.
%
\end{abstract}
\par
\maketitle
{\bf Key words:} liftable vector field, reduced
Kodaira-Spencer-Mather map, higher version of the reduced
Kodaira-Spencer-Mather map,
finite multiplicity, corank at most one.
\tableofcontents
\section{Introduction} 
\label{introduction} \noindent Let $S$ be a finite subset of
$\mathbb{K}^n$, where $\mathbb{K}$ \textcolor{black}{is} the real
field $\mathbb{R}$ or the complex field $\mathbb{C}$ and $n$
\textcolor{black}{is} a positive integer. A map-germ $(\mathbb{K}^n,
S) \to (\mathbb{K}^p,0)$ is called a {\it multigerm}, and it is
called a {\it mono-germ} if $S$ consists of only one point. Let
$C_S$ (resp., $C_0$) be the set of \textcolor{black}{analytic} (that
is, \textcolor{black}{real-analytic}  if $\mathbb{K}=\mathbb{R}$ or
holomorphic if $\mathbb{K}=\mathbb{C}$) multigerms of function
$(\mathbb{K}^n,S)\to \mathbb{K}$ (resp., germs of function
$(\mathbb{K}^p, 0)\to \mathbb{K}$)), and let $m_S$ (resp., $m_0$) be
the subset of $C_S$ (resp., $C_0$) consisting of
\textcolor{black}{analytic}  function-germs $(\mathbb{K}^n,S)\to
(\mathbb{K},0)$ (resp., $(\mathbb{K}^p, 0)\to (\mathbb{K},0)$). It
is clear that the sets $C_S$ and $C_0$ have natural
$\mathbb{K}$-algebra structures induced by the $\mathbb{K}$-algebra
structure of $\mathbb{K}$. For an \textcolor{black}{analytic}
multigerm $f: (\mathbb{K}^n,S)\to (\mathbb{K}^p,0)$, let $f^*:
C_0\to C_S$ be the $\mathbb{K}$-algebra homomorphism defined by
$f^*(u)=u\circ f$.   Set $Q(f)=C_S/f^*m_0C_S$.       A multigerm $f:
(\mathbb{K}^n,S)\to (\mathbb{K}^p,0)$ is said to \textcolor{black}
{have} {\it finite multiplicity} if $Q(f)$ is a finite dimensional
$\mathbb{K}$-vector space. It is well-known that if a multigerm $f:
(\mathbb{K}^n,S)\to (\mathbb{K}^p,0)$ \textcolor{black} {has} finite
multiplicity, then $n$ must be less than or equal to $p$.
\par
For an \textcolor{black}{analytic}  multigerm
$f: (\mathbb{K}^n, S)\to (\mathbb{K}^p,T)$ such that $f(S)\subset
T$, where $S$ (resp., $T$) is a finite subset of $\mathbb{K}^n$
 (resp., $\mathbb{K}^p$),
let $\theta_S(f)$ be the $C_S$-module consisting of germs of
\textcolor{black}{analytic}  vector fields along $f$.
We may identify $\theta_S(f)$ with $\underbrace{C_S\oplus \cdots
\oplus C_S}_{p \mbox{ tuples}}$. We set
$\theta_S(n)=\theta_S(id._{({\mathbb{K}^n},S)})$ and
$\theta_{0}(p)=\theta_{\{0\}}(id._{(\mathbb{K}^p,0)})$, where
$id._{({\mathbb{K}^n},S)}$ (resp., $id._{(\mathbb{K}^p,0)}$) is the
germ of the identity mapping of
$(\mathbb{K}^n,S)$ (resp., $(\mathbb{K}^p,0)$).
\par
For a given \textcolor{black}{analytic}  multigerm $f:
(\mathbb{K}^n,S)\to (\mathbb{K}^p,0)$,
following Mather (\cite{mather3}), we define $tf$ and  $\omega f$ as
follows:
\begin{eqnarray*}
tf: \theta_S(n)\to \theta_S(f),& \quad & tf(\eta)=df\circ \eta, \\
\omega f : \theta_{0}(p) \to \theta_S(f), & \quad & \omega f(\xi)=
\xi\circ f,
\end{eqnarray*}
where $df$ is the differential of $f$. For $f$, following Wall
(\cite{wall}), we set
\begin{eqnarray*}
T\mathcal{R}(f)  =  tf(m_S\theta_S(n)), & { } & T\mathcal{R}_e(f)  =  tf(\theta_S(n)), \\
T\mathcal{L}(f)  =  \omega f(m_0\theta_{0}(p)), & { } & T\mathcal{L}_e(f)  =  \omega f(\theta_{0}(p)), \\
T\mathcal{A}(f)  =  T\mathcal{R}(f)+T\mathcal{L}(f), & { } & T\mathcal{A}_e(f)  =  T\mathcal{R}_e(f)+T\mathcal{L}_e(f), \\
T\mathcal{K}(f)  = T\mathcal{R}(f)+f^*m_0\theta_S(f), & { } &
T\mathcal{K}_e(f)  = T\mathcal{R}_e(f)+f^*m_0\theta_S(f).
\end{eqnarray*}
\par
For a given \textcolor{black}{analytic}  multigerm $f:
(\mathbb{K}^n, S)\to (\mathbb{K}^p,0)$, following Arnol'd
(\cite{arnold}), we call a vector field $\xi\in \theta_0(p)$ {\it
liftable over $f$}  if $\xi\circ f$ belongs to $T\mathcal{R}_e(f)$.
The set of vector fields liftable over $f$ is denoted by $Lift(f)$.
It is clear that $Lift(f)$ naturally has a $C_0$-module structure.
\par
\medskip
\textcolor{black} {The use of liftable vector fields has proven to
be a fundamental tool in the study of classification techniques. In
\cite{arnold} and \cite{damon}, and more recently in
\cite{coopermondwikatique}, \cite{houston} and
\cite{raulcidinharoberta}, $Lift(f)$ has played a central role in
the development of operations and in order to calculate the
codimensions of the multigerms resulting from these operations.}     
\textcolor{black}{However, in general, obtaining generators of $Lift(f)$ can be a very 
hard task. In fact, some articles such as \cite{houstonlittlestone} are devoted to 
constructing $Lift(f)$ for a particular case of germs.}

The purpose of this paper is giving a systematic method to construct
liftable vector fields over a multigerm. In order to create the
systematic method, we first concentrate on obtaining a reasonable
class of multigerms for which the following problems can be
affirmatively answered.
\begin{problem}\label{problem 1}
{\rm Let $f: (\mathbb{K}^n,S)\to (\mathbb{K}^p,0)$ be an
\textcolor{black}{analytic}  multigerm.
\begin{enumerate}
\item Is the module of vector fields liftable over $f$ finitely generated ?
\item Can we characterize the minimal number of generators when the module of vector fields liftable over $f$ is finitely generated ?
\item Can we calculate the minimal number of generators when the module of vector fields liftable over $f$ is finitely generated ?
\item Can we construct generators when the module of vector fields liftable over $f$ is finitely generated ?
\end{enumerate}
}
\end{problem}
\noindent In order to obtain such a reasonable class, we generalize
Mather's homomorphism (\cite{mather4})
\[
\overline{\omega}f: \frac{\theta_0(p)}{m_0\theta_0(p)}\to
\frac{\theta_S(f)}{T\mathcal{K}_e(f)}
\]
defined by $\overline{\omega}f([\xi])=[\omega f(\xi)]$. Notice that
\[
\frac{\theta_S(f)}{T\mathcal{K}_e(f)} \cong
\frac{\frac{\theta_S(f)}{T\mathcal{R}_e(f)}}{f^*m_0\left(\frac{\theta_S(f)}{T\mathcal{R}_e(f)}\right)}
\]
as finite dimensional vector spaces over $\mathbb{K}$ for any
\textcolor{black}{analytic}  multigerm $f$ satisfying
$\dim_{\mathbb{K}}\theta_S(f)/T\mathcal{K}_e(f)<\infty$. Thus, by
the preparation theorem (for instance, see
\cite{arnoldguseinzadevarchenko}), we have that
$\theta_S(f)=T\mathcal{A}_e(f)$ if and only if $\overline{\omega}f$
is surjective for any \textcolor{black}{analytic}  multigerm $f$
satisfying $\dim_{\mathbb{K}}\theta_S(f)/T\mathcal{K}_e(f)<\infty$.
In the case that $\mathbb{K}=\mathbb{C}$, $n\ge p$ and
$S=\{\mbox{one point}\}$, the map $\hat{\omega}f: \theta_0(p)\to
\frac{\theta_S(f)}{T\mathcal{R}_e(f)}$ given by
$\hat{\omega}f(\xi)=[\omega f(\xi)]$ is called the {\it
Kodaira-Spencer map} of $f$ and Mather's homomorphism
$\overline{\omega}f$ is called the {\it reduced Kodaira-Spencer map}
of $f$ (\cite{looijenga}). Thus, $\overline{\omega}f$, which we call
the {\it reduced Kodaira-Spencer-Mather map}, is a generalization of
the reduced Kodaira-Spencer map of $f$; and the module of vector
fields liftable over $f$ is the kernel of $\hat{\omega}f$. We would
like to have higher versions of $\overline{\omega}f$. For a
non-negative integer $i$, an element of $m_S^i$ or $m_0^i$ is a germ
of \textcolor{black}{analytic}  function such that the terms of the
Taylor series of it up to $(i-1)$ are zero. Thus, $m_S^0=C_S$ and
$m_0^0=C_0$. For any non-negative integer $i$ and a given
\textcolor{black}{analytic}  multigerm  $f: (\mathbb{K}^n,S) \to
(\mathbb{K}^p,0)$,
we let
\[
{}_i\overline{\omega}f:
\frac{m_0^i\theta_0(p)}{m_0^{i+1}\theta_0(p)}\to
\frac{f^*m_0^i\theta_S(f)}{T\mathcal{R}_e(f)\cap
f^*m_0^i\theta_S(f)+f^*m_0^{i+1}\theta_S(f)}
\]
be a homomorphism of $C_0$-modules via $f$ defined by
${}_i\overline{\omega}f([\xi])=[\omega f(\xi)]$. Then,
${}_i\overline{\omega}f$ is clearly well-defined. In this paper, we
call ${}_i\overline{\omega}f$ a {\it higher version of reduced
Kodaira-Spencer-Mather map}. Note that
${}_0\overline{\omega}f=\overline{\omega}f$. Similarly as the target
module of $\overline{\omega}f$, for any non-negative integer $i$ and
any \textcolor{black}{analytic}  multigerm $f$ satisfying
$\dim_{\mathbb{K}}\theta_S(f)/T\mathcal{K}_e(f)<\infty$, the target
module of ${}_i\overline{\omega}f$ is isomorphic to the following:
\[
\frac{\frac{f^*m_0^i\theta_S(f)}{T\mathcal{R}_e(f)\cap
f^*m_0^i\theta_S(f)}}
{f^*m_0\left(\frac{f^*m_0^i\theta_S(f)}{T\mathcal{R}_e(f)\cap
f^*m_0^i\theta_S(f)}\right)}.
\]
Thus, again by the preparation theorem, we have that
$f^*m_0^i\theta_S(f)\subset T\mathcal{A}_e(f)$ if and only if
${}_i\overline{\omega}f$ is surjective. The following clearly holds:
\begin{lemma}\label{lemma 1.1}
Let $f: (\mathbb{K}^n,S)\to (\mathbb{K}^p,0)$ be an
\textcolor{black}{analytic}  multigerm satisfying the condition
$\dim_{\mathbb{K}}\theta_S(f)/T\mathcal{K}_e(f)<\infty$. Then, the
following hold:
\begin{enumerate}
\item
Suppose that there exists a non-negative integer $i$ such that
${}_i\overline{\omega}f$ is surjective. Then,
${}_j\overline{\omega}f$ is surjective for any integer $j$ such that
$i< j$.
\item Suppose that there exists a non-negative integer $i$ such that ${}_i\overline{\omega}f$ is injective.
Then, ${}_j\overline{\omega}f$ is injective for any non-negative
integer $j$ such that $i> j$.
\end{enumerate}
\end{lemma}
\begin{definition}\label{definition 1.1}
{\rm Let $f: (\mathbb{K}^n,S)\to (\mathbb{K}^p,0)$ be an
\textcolor{black}{analytic}  multigerm satisfying the condition
$\dim_{\mathbb{K}}\theta_S(f)/T\mathcal{K}_e(f)<\infty$.
\begin{enumerate}
\item Set $
I_1(f)=\left\{i\in \{0\}\cup \mathbb{N}\; |\; {}_i\overline{\omega}f
\mbox{ is surjective}\right\} $. Define $i_1(f)$ as
\[
i_1(f)= \left\{
\begin{array}{cc}
\infty & (\mbox{if }I_1(f)=\emptyset ) \\
\min I_1(f) &  (\mbox{if }I_1(f)\ne \emptyset ).
\end{array}
\right.
\]
\item Set $
I_2(f)=\left\{i\in \{0\}\cup \mathbb{N}\; |\; {}_i\overline{\omega}f
\mbox{ is injective}\right\}. $ Define $i_2(f)$ as
\[
i_2(f)= \left\{
\begin{array}{cc}
-\infty & (\mbox{if }I_2(f)=\emptyset ) \\
\max I_2(f) &  (\mbox{if }\emptyset\ne I_2(f)\ne \{ 0\}\cup\mathbb{N} ) \\
\infty & (\mbox{if }I_2(f)= \{ 0\}\cup\mathbb{N} ).
\end{array}
\right.
\]
\end{enumerate}
}
\end{definition}
\textcolor{black}{An analytic}  multigerm $f: (\mathbb{K}^n,S)\to
(\mathbb{K}^p,0)$ is said to be {\it finitely determined}  if there
exists a positive integer $k$ such that the inclusion
$m_S^k\theta_S(f)\subset T\mathcal{A}_e(f)$ holds. The proof of the
assertion (ii) of proposition 4.5.2 in \cite{wall} works well to
show the following:
\begin{proposition}\label{proposition 1}
Let $f: (\mathbb{K}^n,S)\to (\mathbb{K}^p,0)$ be a finitely
determined multigerm satisfying $\theta_S(f)\ne T\mathcal{A}_e(f)$.
Then, $i_2(f)\ge 0$.
\end{proposition}
From here we concentrate on dealing with the case $n\le p$ because
the purpose of this paper is to construct liftable vector fields
over a multigerm with finite multiplicity. However, the last section
is devoted to extending the results to the case $n>p$.
Suppose that $f: (\mathbb{K}^n,S)\to (\mathbb{K}^p,0)$ is finitely
determined. Then, since it is clear that $f^*m_0C_S\subset m_S$,
there exists a positive integer $k$ such that the inclusion
$f^*m_0^k\theta_S(f)\subset T\mathcal{A}_e(f)$ holds.   Thus,
${}_k\overline{\omega}f$ is surjective.     Conversely, suppose that
there exists a positive integer $k$ such that
${}_k\overline{\omega}f$ is surjective for an
\textcolor{black}{analytic}  multigerm $f: (\mathbb{K}^n,S)\to
(\mathbb{K}^p,0)$ satisfying
$\dim_{\mathbb{K}}\theta_S(f)/T\mathcal{K}_e(f)<\infty$. Then, as we
have already confirmed, the inclusion $f^*m_0^k\theta_S(f)\subset
T\mathcal{A}_e(f)$ holds by the preparation theorem. In the case
$n\le p$, by Wall's estimate (theorem 4.6.2 in \cite{wall}), the
condition $\dim_{\mathbb{K}}\theta_S(f)/T\mathcal{K}_e(f)<\infty$
implies that there exists an integer $\ell$ such that
$m_S^\ell\subset f^*m_0C_S$.     Hence, we have the following:
\begin{proposition}\label{proposition 2}
Let $f: (\mathbb{K}^n,S)\to (\mathbb{K}^p,0)$ be an
\textcolor{black}{analytic}  multigerm satisfying the condition
$\dim_{\mathbb{K}}\theta_S(f)/T\mathcal{K}_e(f)<\infty$. Suppose
that $n\le p$. Then, $i_1(f)<\infty$ if and only if $f$ is finitely
determined.
\end{proposition}
\textcolor{black}{An analytic}  multigerm $f: (\mathbb{K}^n, S)\to
(\mathbb{K}^p,0)$ $(n\le p)$ is said to be {\it of corank at most
one} if $\max\{n-\mbox{rank} Jf(s_j)\; |\; 1\le j\le |S|\}\le 1 $
holds, where $Jf(s_j)$ is the Jacobian matrix of $f$ at $s_j\in S$
and $|S|$ stands for the number of distinct points of $S$.
\begin{proposition}\label{proposition 3}
Let $f: (\mathbb{K}^n,S)\to (\mathbb{K}^p,0)$ $(n\le p)$ be a
finitely determined multigerm of corank at most one.
Then, $i_1(f)\ge i_2(f)$.
\end{proposition}
\noindent Proposition \ref{proposition 3} is proved in \S
\ref{section 2}.     Proposition \ref{proposition 3} yields the
following corollary.
\begin{corollary}\label{corollary 1}
Let $f: (\mathbb{K}^n,S)\to (\mathbb{K}^p,0)$ $(n\le p)$ be a
finitely determined multigerm of corank at most one. Suppose that
there exists a non-negative integer $i$ such that $i_1(f)=i_2(f)=i$.
Then, the following hold:
\begin{enumerate}
\item For any non-negative integer $j$ such that $j< i$,
${}_j\overline{\omega}f$ is injective but not surjective.
\item For any non-negative integer $j$ such that
$i<j$, ${}_j\overline{\omega}f$ is surjective but not injective.
\end{enumerate}
\end{corollary}
\begin{example}\label{example 1.1}
{\rm Let $e: \mathbb{K}\to \mathbb{K}^2$ be the embedding defined by
$e(x)=(x, 0)$ and for any real number $\theta$ let $R_{\theta}:
\mathbb{K}^2\to \mathbb{K}^2$ be the linear map which gives the
rotation of $\mathbb{K}^2$ about the origin with respect to the
angle $\theta$.
$$
R_{\theta} \left(
\begin{array}{c}
X \\
Y
\end{array}
\right) = \left(
\begin{array}{rr}
\cos\theta & -\sin\theta \\
\sin\theta & \cos\theta
\end{array}
\right) \left(
\begin{array}{c}
X \\
Y
\end{array}
\right).
$$
For any non-negative integer $\ell$ set $S=\{s_0, \ldots,
s_{\ell+1}\}$ $(s_j\ne s_k \mbox{ if }j\ne k)$. Define
$\theta_j=j\frac{\pi}{\ell+2}$ and set $e_{j}: (\mathbb{K}, s_j)\to
(\mathbb{K}^2,0)$ as $e_{j}(x_j)=R_{\theta_j}\circ e(x_j)$ for any
$j$ $(0\le j\le \ell+1)$, where $x_j=x-s_j$.     Then,
$E_\ell=\{e_0, \ldots, e_{\ell+1}\}: (\mathbb{K},S)\to
(\mathbb{K}^2,0)$ is a finitely determined multigerm of corank at
most one. The image of $E_\ell$ is a line arrangement and hence the
Euler vector field of the defining equation of the image of $E_\ell$
is a liftable vector field over $E_\ell$. It follows that
${}_1\overline{\omega}E_\ell$ is not injective. Furthermore, it is
easily seen that ${}_0\overline{\omega}E_\ell$ is injective even in
the case $\ell=0$ (in the case $\ell\ge 1$ this is a trivial
corollary of Proposition \ref{proposition 1}). Thus,
$i_2(E_\ell)=0$. On the other hand, it is not hard to show that
$i_1(E_\ell)=\ell$.       Therefore, $i_1(E_\ell)-i_2(E_\ell)=\ell$.
}
\end{example}
This example shows that, in general, there are no upper bounds of
$i_1(f)-i_2(f)$ for a finitely determined multigerm $f:
(\mathbb{K}^n,S)\to (\mathbb{K}^p,0)$ $(n\le p)$ of corank at most
one. This example shows also that the integer $i_1(f)-i_2(f)$
measures how well-behaved a given finitely determined multigerm of
corank at most one is from the viewpoint of liftable vector fields.
\par
The following Theorem \ref{theorem 1} shows that the desired
reasonable class is the set consisting of finitely determined
multigerms $f: (\mathbb{K}^n, S)\to (\mathbb{K}^p,0)$ $(n\le p)$ of
corank at most one satisfying $i_1(f)=i_2(f)$.
\begin{theorem}\label{theorem 1}
Let $f: (\mathbb{K}^n, S)\to (\mathbb{K}^p,0)$ $(n\le p)$ be a
finitely determined multigerm of corank at most one. Suppose that
there exists a non-negative integer $i$ such that $i_1(f)=i_2(f)=i$.
Then,
the minimal number of generators for the module of vector fields
liftable over $f$ is exactly $\dim_{\mathbb{K}}\mbox{\rm
ker}({}_{i+1}\overline{\omega}f)$.
\end{theorem}
\noindent Notice that the embedding $e$ in Example \ref{example 1.1}
does not satisfy the assumption of Theorem \ref{theorem 1}.
Actually, since ${}_0\overline{\omega}e$ is surjective but not
injective, $i_1(e)=0$ and $i_2(e)=-\infty$.
On the other hand, the multigerm $E_0$ 
in Example \ref{example 1.1} does satisfy the assumption of Theorem
\ref{theorem 1} though $E_0$ does not satisfy the assumption of
Proposition \ref{proposition 1}. Furthermore, a lot of examples of
Theorem \ref{theorem 1} are given by Proposition \ref{proposition
new 1} (see also Section \ref{section 3}).
\begin{definition}
\begin{enumerate}
\item A multigerm $f: (\mathbb{K}^n, S)\to (\mathbb{K}^p,0)$ is said to be {\it stable} if
it satisfies $\theta_S(f)=T_e\mathcal{A}(f)$.
\item Define the mapping $ev_0: \theta_0(p)\to T_0(\mathbb{R}^p)$ by
$ev_0(\eta)=\eta(0)$.
\item A stable multigerm $f: (\mathbb{K}^n, S)\to (\mathbb{K}^p,0)$ is said to be {\it isolated}
if $ev_0(\eta)=0$ for any $\eta\in Lift(f)$.
\end{enumerate}
\end{definition}
\noindent The following proposition shows that our reasonable class
contains the set consisting of isolated stable multigerms $f:
(\mathbb{K}^n, S)\to (\mathbb{K}^p,0)$ $(n\le p)$ of corank at most
one.
\begin{proposition}\label{proposition new 1}
Let $f: (\mathbb{K}^n, S)\to (\mathbb{K}^p,0)$ $(n\le p)$ be a
finitely determined multigerm of corank at most one. Then, the
following hold:
\begin{enumerate}
\item In the case $i=0$, the following hold:
\begin{enumerate}
\item
${}_0\overline{\omega}f$ is surjective if and only if $f$ is stable.
\item
${}_0\overline{\omega}f$ is injective if and only if $f$ is
isolated.
\end{enumerate}
\item In the case $i=1$, the following hold:
\begin{enumerate}
\item
${}_1\overline{\omega}f$ is surjective if and only if
$T\mathcal{A}(f)=T\mathcal{K}(f)$.
\item
${}_1\overline{\omega}f$ is injective if only if for any $\eta\in
Lift(f)$ $\eta$ has no constant terms and no linear terms.
Moreover, these equivalent conditions imply that
$\dim_{\mathbb{K}}\theta_S(f)/T\mathcal{A}_e(f)>1$.
\end{enumerate}
\end{enumerate}
\end{proposition}
\par
\medskip
Next, in order to answer (3) of Problem \ref{problem 1} for a given
finitely determined multigerm $f: (\mathbb{K}^n, S)\to
(\mathbb{K}^p,0)$ $(n\le p)$ of corank at most one such that $0\le
i_1(f)=i_2(f)<\infty$,
we generalize Wall's homomorphism (\cite{wall})
\[
\overline{t} f : Q(f)^n\to Q(f)^p,\quad \overline{t}f([\eta
])=[tf(\eta )]
\]
as follows. For a given \textcolor{black}{analytic}  multigerm $f:
(\mathbb{K}^n, S)\to (\mathbb{K}^p,0)$ satisfying the condition
$\dim_{\mathbb{K}}Q(f)<\infty$, let $\delta(f)$ (resp., $\gamma(f)$)
be the dimension of the vector space $Q(f)$ (resp., the dimension of
the kernel of $\overline{t}f$).   For the $f$
and a non-negative integer $i$, we Set ${ }_iQ(f) =
f^*m_0^iC_S/f^*m_0^{i+1}C_S$ and
${}_i\delta(f)=\dim_{\mathbb{K}}{}_iQ(f)$. Thus, we have that
${}_0Q(f)=Q(f)$ and ${}_0\delta(f)=\delta(f)=\dim_{\mathbb{K}}Q(f)$.
The $Q(f)$-modules ${}_iQ(f)^n$ and ${}_iQ(f)^p$ may be identified
with the following respectively.
\[
\frac{f^*m_0^i\theta_S(n)}{f^*m_0^{i+1}\theta_S(n)}
\quad\mbox{and}\quad
\frac{f^*m_0^i\theta_S(f)}{f^*m_0^{i+1}\theta_S(f)}.
\]
Let ${}_i\gamma(f)$ be the dimension of the kernel of the following
well-defined homomorphism of $Q(f)$-modules.
\[
{}_i\overline{t}f : {}_iQ(f)^n \to {}_iQ(f)^p,\quad
{}_i\overline{t}f([\eta ])=[tf(\eta )].
\]
Then, we have that ${}_i\delta(f)<\infty$ if $\delta(f) < \infty$
and ${}_i\gamma(f)<\infty$ if $\gamma(f) < \infty$. For details on
${}_iQ(f)$, ${}_i\delta(f)$, ${}_i\overline{t}f$ and
${}_i\gamma(f)$, see \cite{nishimura}.
\begin{proposition}\label{proposition 4}
Let $f: (\mathbb{K}^n, S)\to (\mathbb{K}^p,0)$ be an
\textcolor{black}{analytic}  multigerm with finite multiplicity and
of corank at most one.
Suppose that there exists a non-negative integer
$i$ such that ${}_{i+1}\overline{\omega}f$ is surjective. Then, the
following holds:
\[
\dim_{\mathbb{K}}\mbox{\rm ker}({}_{i+1}\overline{\omega}f)= p\cdot
\left(
\begin{array}{c}
p+i \\
i+1
\end{array}
\right)
-\left((p-n)\cdot
{}_{i+1}\delta(f)+{}_{i+1}\gamma(f)-{}_i\gamma(f)\right),
\]
where the dot in the center stands for the multiplication.
\end{proposition}
\begin{proposition}\label{proposition 5}
Let $f: (\mathbb{K}^n, S)\to (\mathbb{K}^p,0)$ be an
\textcolor{black}{analytic}  multigerm with finite multiplicity and
of corank at most one. Then, the following hold:
\begin{enumerate}
\item ${}_0\gamma(f)=\gamma(f)=\delta(f)-|S|$.
\item
\[
{}_i\delta(f)= \left(
\begin{array}{c}
n+i-1 \\
i
\end{array}
\right) \cdot \delta(f),\; {}_i\gamma(f)= \left(
\begin{array}{c}
n+i-1 \\
i
\end{array}
\right) \cdot \gamma(f) \quad (i\in \mathbb{N}\cup \{0\}).
\]
\end{enumerate}
\end{proposition}
By combining Propositions \ref{proposition 4} and \ref{proposition
5}, for an \textcolor{black}{analytic}  multigerm $f$ of corank at
most one such that $\dim_{\mathbb{K}}Q(f)<\infty$, the
$\mathcal{A}$-invariant \lq\lq $\dim_{\mathbb{K}}\mbox{\rm
ker}({}_{i+1}\overline{\omega}f)$\rq\rq can be calculated easily by
using $\mathcal{K}$-invariants \lq\lq $\delta(f), \gamma(f)$\rq\rq
when there exists a non-negative integer $i$ such that
${}_{i+1}\overline{\omega}f$ is surjective.
\begin{theorem}\label{theorem 2}
Let $f: (\mathbb{K}^n, S)\to (\mathbb{K}^p,0)$ be an
\textcolor{black}{analytic}  multigerm with finite multiplicity and
of corank at most one.
\begin{enumerate}
\item
Let $F: (\mathbb{K}^n\times \mathbb{K}^r, S\times \{0\})\to
(\mathbb{K}^p\times \mathbb{K}^r,(0,0))$ be a stable unfolding of
$f$. Let $\eta=(\eta_1, \ldots, \eta_p, \eta_{p+1},
\textcolor{black}{\ldots, }\eta_{p+r})$ be an element of the
intersection $Lift(F)\cap Lift(g)$, where $g=\{g_1, \ldots, g_r\}$
with $g_i: (\mathbb{K}^p\times \mathbb{K}^r, (0,0))\to
(\mathbb{K}^p\times \mathbb{K}^r,(0,0))$ defined by $g_i(X_1,
\ldots, X_p, \lambda_1, \ldots, \lambda_r)= (X_1, \ldots, X_p,
\lambda_1, \ldots, \lambda_{i-1}, $ $\lambda_i^2, \lambda_{i+1},
\ldots, \lambda_r)$ $(1\le i\le r)$. Then,
$\overline{\eta}(X)=(\eta_1(X, 0), \ldots, \eta_p(X,0))$ is an
element of $Lift(f)$.
\item
Suppose that $f$ admits a one-parameter stable unfolding $F:
(\mathbb{K}^n\times \mathbb{K}, S\times \{0\})\to
(\mathbb{K}^p\times \mathbb{K},(0,0))$. Then, for any
$\overline{\eta}\in Lift(f)$ there exists an element $\eta=(\eta_1,
\ldots, \eta_p, \eta_{p+1}) \in Lift(F)\cap Lift(g_1)$ such that the
equality $\textcolor{black} {\overline{\eta}}(X)=(\eta_1(X,0),
\ldots, \eta_p(X,0))$ holds.
\end{enumerate}
\end{theorem}
\noindent
The proof of Theorem \ref{theorem 1} provides a recipe for
constructing all liftable vector fields over a finitely determined
multigerm $f$ of corank at most one satisfying $i_1(f)=i_2(f)$. In
particular, by Proposition \ref{proposition new 1}, all liftable
vector fields over an isolated stable multigerm $f$ of corank at
most one can be constructed. \textcolor{black} {Since any stable
germ is $\mathcal A$-equivalent to a prism on an isolated stable
multigerm, all liftable vector fields over any stable germ can be
constructed.} Moreover, by using Theorem \ref{theorem 2}, we can
construct all liftable vector fields over an
\textcolor{black}{analytic} multigerm $f$ of corank at most one
admitting a one-parameter stable unfolding $F$ from $Lift(F)$. It is
clear also that if $f$ satisfies
$\dim_{\mathbb{K}}\theta_S(f)/T\mathcal{A}_e(f)= 1$ (namely, $f$ is
a multigerm of  $\mathcal{A}_e$-codimension one), then $f$ admits a
one-parameter stable unfolding. Thus, we can construct all liftable
vector fields over a multigerm of corank at most one and of
$\mathcal{A}_e$-codimension one. In particular, for any augmentation
defined in \cite{coopermondwikatique}, all liftable vector fields
over it can be constructed by our recipe. In Remark \ref{remopsu} at
the end of Section \ref{section 6} an idea on how big the space of
germs which admit a one-parameter stable unfolding is is given.
\par
\smallskip
This paper is organized as follows. In Section \ref{section 2},
proofs of Propositions \ref{proposition 3}, \ref{proposition new 1},
\ref{proposition 4}, and \ref{proposition 5} are given. In Section
\ref{section 3}, examples for which actual calculations of minimal
numbers of generators are carried out
are given. Theorem \ref{theorem 1} (resp., Theorem \ref{theorem 2})
is proved in Sections \ref{section 4} (resp., Section \ref{section
5}). In Section \ref{section 6}, by constructing concrete generators
for several examples using Theorem \ref{theorem 2} and the proof of
Theorem \ref{theorem 1}, it is explained in detail how to construct
liftable vector fields over an \textcolor{black}{analytic} multigerm
of corank at most one admitting a one-parameter stable unfolding.
Finally, Section \ref{section 7} generalizes the results for the
case $n>p$.
\section{Proofs of Propositions \ref{proposition 3}, \ref{proposition new 1},
\ref{proposition 4}, and \ref{proposition 5}} \label{section 2}
Firstly, Proposition \ref{proposition 5} is proved.
\par
\smallskip
\noindent \underline{\it Proof of Proposition \ref{proposition
5}.}\qquad
\par
Set $S=\{s_1, \ldots, s_{|S|}\}$ $(s_j\ne s_k \mbox{ if }j\ne k)$
and for any $j$ $(1\le j\le |S|)$ let $f_j$ be the restriction
$f|_{(\mathbb{K}^n,s_j)}$. Then, we have the following:
\begin{eqnarray*}
\delta(f)=\dim_{\mathbb{K}}Q(f)=\sum_{j=1}^{|S|}\dim_{\mathbb{K}}Q(f_j)
& = & \sum_{j=1}^{|S|}\delta(f_j).  \\
\gamma(f)=\dim_{\mathbb{K}}\mbox{\rm ker}(\overline{t}f)
=\sum_{j=1}^{|S|}\dim_{\mathbb{K}}\mbox{\rm ker}(\overline{t}f_j)
& = & \sum_{j=1}^{|S|}\gamma(f_j) \\
{ } & = & \sum_{j=1}^{|S|}\left(\delta(f_j)-1\right)=\delta(f)-|S|.
\end{eqnarray*}
This completes the proof of the assertion 1 of Proposition
\ref{proposition 5}.
\par
Next we prove the assertion 2 of Proposition \ref{proposition 5}.
Since $f$ is of corank at most one, for any $j$ $(1\le j\le |S|)$
there exist germs of diffeomorphism $h_j: (\mathbb{K}^n,s_j)\to
(\mathbb{K}^n,s_j)$ and  $H_j: (\mathbb{K}^p,0)\to (\mathbb{K}^p,0)$
such that $H_j\circ f_j\circ h_j^{-1}$ has the following form:
\[
H_j\circ f_j\circ h_j^{-1}(x, y)
 =  (x, y^{\delta(f_j)}+f_{j, n}(x, y), f_{j, n+1}(x, y), \ldots, f_{j, p}(x, y)).
\]
Here, $(x,y)=(x_1, \ldots, x_{n-1}, y)$ is the local coordinate with
respect to the coordinate neighborhood $(U_j, h_j)$ at $s_j$ and
$f_{j, q}$ satisfies $f_{j,q}(0, \ldots, 0, y)=o(y^{\delta(f_j)})$
for any $q$ $(n\le q\le p)$.
By the preparation theorem, $C_{s_j}$ is generated by $1, y, \ldots,
y^{\delta(f_j)-1}$ as $C_0$-module via $f_j$. Thus,
$f_j^*m_0^iC_{s_j}$ is generated by elements of the following set as
$C_0$-module via $f_j$.
\[
\left\{x_1^{k_1}\cdots x_{n-1}^{k_{n-1}}y^{k_n\delta(f_j)+\ell}\;
\left|\; k_m\ge 0, \sum_{m=1}^nk_m=i, 0\le \ell\le
\delta(f_j)-1\right.\right\}.
\]
Thus, the following set is a basis of ${}_iQ(f_j)$.
\[
\left.\left\{\left[ x_1^{k_1}\cdots
x_{n-1}^{k_{n-1}}y^{k_n\delta(f_j)+\ell}\right]\; \left|\; k_m\ge 0,
\sum_{m=1}^nk_m=i, 0\le \ell\le \delta(f_j)-1\right.\right.\right\}.
\]
Therefore, we have the following:
\begin{eqnarray*}
{}_i\delta(f)=\dim_{\mathbb{K}}{}_iQ(f)
=\sum_{j=1}^{|S|}\dim_{\mathbb{K}}{}_iQ(f_j)
& = & \sum_{j=1}^{|S|}\dim_{\mathbb{K}}\frac{f_j^*m_0^iC_{s_j}}{f_j^*m_0^{i+1}C_{s_j}} \\
{ } & = & \sum_{j=1}^{|S|} \left(
\begin{array}{c}
n+i-1 \\
i
\end{array}
\right) \cdot
\delta(f_j) \\
{ } & = & \left(
\begin{array}{c}
n+i-1 \\
i
\end{array}
\right) \cdot
\sum_{j=1}^{|S|}\delta(f_j) \\
{ } & = & \left(
\begin{array}{c}
n+i-1 \\
i
\end{array}
\right) \cdot \delta(f).
\end{eqnarray*}
\par
\smallskip
Next we prove the formula for ${}_i\gamma(f)$.    Since it is clear
that ${}_i\gamma(f_j)$ does not depend on the particular choice of
coordinate systems of $(\mathbb{K}^n,s_j)$ and of
$(\mathbb{K}^p,0)$, we may assume that $f_j$ has the above form from
the first. Then, it is easily seen that the following set is a basis
of $\mbox{\rm ker}{}_i\overline{t}f_j$. {\small
\[
\left\{\left[\underbrace{0\oplus\cdots\oplus 0}_{(n-1)\mbox{
tuples}}\oplus x_1^{k_1}\cdots
x_{n-1}^{k_{n-1}}y^{k_n\delta(f_j)+\ell}\right]\; \left|\; k_m\ge 0,
\sum_{m=1}^nk_m=i, 1\le \ell\le \delta(f_j)-1\right.\right\}.
\]
} Therefore, we have the following:
\begin{eqnarray*}
{}_i\gamma(f)=\dim_{\mathbb{K}}\mbox{\rm ker}({}_i\overline{t}f)
=\sum_{j=1}^{|S|}\dim_{\mathbb{K}}\mbox{\rm
ker}({}_i\overline{t}f_j)
{ } & = & \sum_{j=1}^{|S|} \left(
\begin{array}{c}
n+i-1 \\
i
\end{array}
\right) \cdot
\left(\delta(f_j)-1\right) \\
{ } & = & \left(
\begin{array}{c}
n+i-1 \\
i
\end{array}
\right) \cdot
\sum_{j=1}^{|S|}\gamma(f_j) \\
{ } & = & \left(
\begin{array}{c}
n+i-1 \\
i
\end{array}
\right) \cdot \gamma(f).
\end{eqnarray*}
\hfill Q.E.D.
\par
\medskip
Secondly, Proposition \ref{proposition 4} is proved.
\par
\smallskip
\noindent \underline{\it Proof of Proposition \ref{proposition
4}}\qquad
\par
Consider the linear map ${}_{i+1}\overline{t}f$. Then,
we have the following:
\[
\dim_{\mathbb{K}}{}_{i+1}Q(f)^n={}_{i+1}\gamma(f)+\dim_{\mathbb{K}}\mbox{\rm
\mbox{\rm Image}}({}_{i+1}\overline{t}f).
\]
Since $\dim_{\mathbb{K}}Q(f)<\infty$ and $f$ is of corank at most
one, it is easily seen that $tf$ is injective.   Hence we see that
$$
\dim_{\mathbb{K}} \frac{T\mathcal{R}_e(f)\cap
f^*m_0^{i+1}\theta_S(f)}{T\mathcal{R}_e(f)\cap
f^*m_0^{i+2}\theta_S(f)} =\dim_{\mathbb{K}}\mbox{\rm \mbox{\rm
Image}}({}_{i+1}\overline{t}f)+{}_i\gamma(f).
$$
Therefore, we have the following:
$$
\dim_{\mathbb{K}}\frac{f^*m_0^{i+1}\theta_S(f)}{T\mathcal{R}_e(f)\cap
f^*m_0^{i+1}\theta_S(f)+f^*m_0^{i+2}\theta_S(f)} = (p-n)\cdot
{}_{i+1}\delta(f)+{}_{i+1}\gamma(f)-{}_i\gamma(f).
$$
Hence,
we have the following:
\begin{eqnarray*}
{ } & { } & \dim_{\mathbb{K}}\mbox{\rm ker}\left({}_{i+1}\overline{\omega}f\right) \\
{ } & = &
\dim_{\mathbb{K}}\frac{m_0^{i+1}\theta_0(p)}{m_0^{i+2}\theta_0(p)} -
\dim_{\mathbb{K}}\frac{f^*m_0^{i+1}\theta_S(f)}{T\mathcal{R}_e(f)\cap
f^*m_0^{i+1}\theta_S(f)+f^*m_0^{i+2}\theta_S(f)}
\\
{ } & = & p\cdot \left(
\begin{array}{c}
p+i \\
i+1
\end{array}
\right) - \left((p-n)\cdot
{}_{i+1}\delta(f)+{}_{i+1}\gamma(f)-{}_i\gamma(f) \right).
\end{eqnarray*}
\hfill Q.E.D.
\par
\medskip
Thirdly, Proposition \ref{proposition 3} is proved.
\par
\smallskip
\noindent \underline{\it Proof of Proposition \ref{proposition
3}}\qquad
\par
By Lemma \ref{lemma 1.1}, it sufficies to show that for any $i$ and
any finitely determined multigerm $f: (\mathbb{K}^n, S)\to
(\mathbb{K}^p,0)$ $(n\le p)$ of corank at most one satisfying that
${}_i\overline{\omega}f$ is surjective, ${}_{i+1}\overline{\omega}f$
is not injective. By Lemma \ref{lemma 1.1},
${}_{i+1}\overline{\omega}f$ is not injective if
${}_{i}\overline{\omega}f$ is not injective.      Thus, we may
assume that ${}_i\overline{\omega}f$ is bijective.
\par
We first prove Proposition \ref{proposition 3} in the case $i=0$.
Since we have assumed that ${}_0\overline{\omega}f$ is bijective,
the following holds (see \cite{mather6} or \cite{wall}):
\[
p\cdot \left(
\begin{array}{c}
p-1 \\
0
\end{array}
\right) = \dim_{\mathbb{K}}\frac{\theta_0(p)}{m_0\theta_0(p)} =
(p-n)\cdot {}_0\delta(f)+{}_0\gamma(f).
\]
Note that the above equality can not be obtained by Proposition
\ref{proposition 4}. Note further that at least one of $p-n>0$ or
${}_0\gamma(f)>0$ holds by this equality. We have the following:
\begin{eqnarray*}
p\cdot \left(
\begin{array}{c}
p \\
1
\end{array}
\right) & = & p^2\cdot \left(
\begin{array}{c}
p-1 \\
0
\end{array}
\right) \\
{ } & = &
p\cdot \left((p-n)\cdot {}_0\delta(f)+{}_0\gamma(f)\right) \\
{ } & = &
\frac{p}{n}\cdot \left((p-n)\cdot {}_1\delta(f)+{}_1\gamma(f)\right) \qquad  (\mbox{\rm by 2 of Proposition 5})\\
{ } & \ge &
(p-n)\cdot {}_1\delta(f)+{}_1\gamma(f)  \qquad    (\mbox{by }n\le p)\\
{ } & \ge &
(p-n)\cdot {}_1\delta(f)+{}_1\gamma(f)-{}_0\gamma(f) \qquad (\mbox{by }{}_0\gamma(f)\ge 0).   \\
\end{eqnarray*}
Since we have confirmed that at least one of $p-n>0$ or
${}_0\gamma(f)>0$ holds, we have the following sharp inequality:
\[
p\cdot \left(
\begin{array}{c}
p \\
1
\end{array}
\right)
>
(p-n)\cdot {}_1\delta(f)+{}_1\gamma(f)-{}_0\gamma(f).
\]
Hence ${}_1\overline{\omega}f$ is not injective by Lemma \ref{lemma
1.1} and Proposition \ref{proposition 4}.
\par
Next we prove Proposition \ref{proposition 3} in the case $i\ge 1$.
Since we have assumed that ${}_i\overline{\omega}f$ is bijective, we
have the following equality by Proposition \ref{proposition 4}:
\begin{eqnarray*}
{ } & { } & p\cdot \left(
\begin{array}{c}
p+i-1 \\
i
\end{array}
\right) \\
{ } & = &
(p-n)\cdot {}_i\delta(f)+{}_i\gamma(f)-{}_{i-1}\gamma(f) \\
{ } & = & (p-n)\cdot
{}_i\delta(f)+\left(1-\frac{i}{n+i-1}\right)\cdot {}_i\gamma(f)
\qquad
(\mbox{\rm by 2 of Proposition \ref{proposition 5}})\\
{ } & = &
(p-n)\cdot {}_i\delta(f)+\frac{n-1}{n+i-1}\cdot {}_i\gamma(f). \\
\end{eqnarray*}
Note that at least one of $p-n>0$ or $(n-1)\cdot {}_i\gamma(f)> 0$
holds by this equality. We have the following:
\begin{eqnarray*}
{ } & { } & p\cdot \left(
\begin{array}{c}
p+i \\
i+1
\end{array}
\right) \\
{ } & = & \frac{p+i}{i+1} \cdot p\cdot \left(
\begin{array}{c}
p+i-1 \\
i
\end{array}
\right) \\
{ } & = &
\frac{p+i}{i+1}\cdot \left((p-n)\cdot {}_i\delta(f)+\frac{n-1}{n+i-1}\cdot {}_i\gamma(f)\right) \\
{ } & = & \frac{p+i}{i+1}\cdot \left(\frac{i+1}{n+i}\cdot (p-n)\cdot
{}_{i+1}\delta(f)
+\frac{n-1}{n+i-1}\cdot \frac{i+1}{n+i}\cdot {}_{i+1}\gamma(f)\right) \\
{ } & { } & \qquad\qquad\qquad\qquad\qquad\qquad\qquad
\qquad  (\mbox{\rm by 2 of Proposition 5})\\
{ } & = & \frac{p+i}{n+i}\cdot (p-n)\cdot {}_{i+1}\delta(f)
+\frac{p+i}{n+i-1}\cdot \frac{n-1}{n+i}\cdot {}_{i+1}\gamma(f) \\
{ } & \ge & (p-n)\cdot {}_{i+1}\delta(f)
+\frac{p+i}{n+i-1}\cdot \frac{n-1}{n+i}\cdot {}_{i+1}\gamma(f)  \qquad    (\mbox{by }n\le p) \\
{ } & \ge & (p-n)\cdot {}_{i+1}\delta(f)
+\frac{n-1}{n+i}\cdot {}_{i+1}\gamma(f)  \qquad    (\mbox{by }n\le p \mbox{ and }(n-1){}_{i+1}\gamma(f)\ge 0) \\
{ } & = & (p-n)\cdot {}_{i+1}\delta(f)
+{}_{i+1}\gamma(f)-{}_i\gamma(f) \qquad  (\mbox{\rm by 2 of
Proposition \ref{proposition 5}}).
\end{eqnarray*}
Since we have confirmed that at least one of $p-n>0$ or $(n-1)\cdot
{}_i\gamma(f)> 0$ holds and ${}_{i+1}\gamma(f)=\frac{n+i}{i+1}\cdot
{}_i\gamma(f)$ by the assertion 2 of Proposition \ref{proposition
5}, we have the following sharp inequality:
\[
p\cdot \left(
\begin{array}{c}
p+i \\
i+1
\end{array}
\right)
>
(p-n)\cdot {}_{i+1}\delta(f)+{}_{i+1}\gamma(f)-{}_i\gamma(f).
\]
Hence, ${}_{i+1}\overline{\omega}f$ is not injective by Lemma 1.1
and Proposition \ref{proposition 4}. \hfill Q.E.D.
\par
\medskip
Finally, Proposition \ref{proposition new 1} is proved.
\par
\smallskip
\noindent \underline{\it Proof of Proposition \ref{proposition new
1}}\qquad
\par
Proof of the assertion (1) of Proposition \ref{proposition new 1} is
as follows. Recall that ${}_0\overline{\omega}f$ is nothing but
Mather's $\overline{\omega}f$ defined in \cite{mather4}.      The
assertion (a) has been already shown by Mather (see p.228 of
\cite{mather4}). Since $\overline{\omega}f:
\theta_0(p)/m_0\theta_0(p)\to \theta_S(f)/T\mathcal{K}_e(f)$ is
defined by $\overline{\omega}f([\eta])=[\eta\circ f]$, by
definition, the injectivity of $\overline{\omega}f$ is equivalent to
assert that $ev_0(\eta)=\eta(0)=0$ for any $\eta\in Lift(f)$.
\par
\smallskip
Proof of the assertion (2) of Proposition \ref{proposition new 1} is
as follows. Recall that ${}_1\overline{\omega}f$ is the following
mapping defined by ${}_1\overline{\omega}f([\eta])=[\eta \circ f]$:
\[
{}_1\overline{\omega}f: \frac{m_0\theta_0(p)}{m_0^{2}\theta_0(p)}\to
\frac{f^*m_0\theta_S(f)}{T\mathcal{R}_e(f)\cap
f^*m_0\theta_S(f)+f^*m_0^{2}\theta_S(f)}.
\]
As we have already confirmed in Section \ref{introduction}, by the
preparation theorem, we have that ${}_1\overline{\omega}f$ is
surjective if and only if $f^*m_0\theta_S(f)\subset
T\mathcal{A}_e(f)$.      Since $f$ is finitely determined, we can
conclude that $f^*m_0\theta_S(f)\subset T\mathcal{A}_e(f)$ if and
only if $f^*m_0\theta_S(f)\subset T\mathcal{A}(f)$.
\par
For  (b) of (2), it is easily seen that injectivity of
${}_1\overline{\omega}f$ is equivalent to assert that $\eta$ has no
constant terms and no linear terms for any $\eta\in Lift(f)$. These
equivalent conditions imply that $f$ is not $\mathcal{A}$-equivalent
to a quasi-homogeneous multigerm.     By \cite{coopermondwikatique},
this implies that $\dim_{\mathbb{K}}\theta_S(f)/T\mathcal{A}_e(f) >
1$. \hfill Q.E.D.
\section{Examples of Theorem \ref{theorem 1}}
\label{section 3}
\begin{example}\label{example 3.1}
{\rm Let $\varphi: (\mathbb{K}^n,0)\to (\mathbb{K}^n,0)$ be the
map-germ given by $\varphi(x_1, \ldots, x_{n-1},$ $y)=(x_1, \ldots,
x_{n-1}, y^{n+1}+\sum_{i=1}^{n-1}x_iy^i)$. Then, it is known that
$f$ is an isolated stable mono-germ by \cite{morin} or
\cite{mather5}. Thus, by Proposition \ref{proposition new 1},
${}_0\overline{\omega}\varphi$ is bijective.
Therefore, it follows that $i_1(\varphi)=i_2(\varphi)=0$. By Theorem
\ref{theorem 1}, Lemma \ref{lemma 1.1} and Propositions
\ref{proposition 4}, \ref{proposition 5}, the minimal number of
generators for $Lift(\varphi)$
can be calculated as follows:
\begin{eqnarray*}
{ } & { } & n\cdot \left(
\begin{array}{c}
n \\
1
\end{array}
\right)
-
\left((n-n)\cdot {}_{1}\delta(\varphi)+{}_{1}\gamma(\varphi)-{}_0\gamma(\varphi)\right)  \\
{ } & = &
n^2-\left((n-n)\cdot n\cdot (n+1)+n\cdot (n+1-1)-(n+1-1)\right)  \\
{ } & = & n.
\end{eqnarray*}
It has been verified in \cite{arnold} that the minimal number of
generators for $Lift(\varphi)$
is exactly $n$ in the complex case. }
\end{example}
\begin{example}\label{example 3.2}
{\rm Let $\varphi_k: (\mathbb{K}^{2k-2},0)\to (\mathbb{K}^{2k-1},0)$
be given by
\begin{eqnarray*}
{ } & { } & \varphi_k(u_1, \ldots, u_{k-2}, v_1, \ldots, v_{k-1}, y) \\
{ } & = & \left(u_1, \ldots, u_{k-2}, v_1, \ldots, v_{k-1},
y^k+\sum_{i=1}^{k-2}u_iy^i, \sum_{i=1}^{k-1}v_iy^i\right).
\end{eqnarray*}
Then, it is known that $f$ is an isolated stable mono-germ by
\cite{morin} or \cite{mather5}. Thus, by Proposition
\ref{proposition new 1}, ${}_0\overline{\omega}\varphi_k$ is
bijective. Therefore, it follows that
$i_1(\varphi_k)=i_2(\varphi_k)=0$. By Theorem \ref{theorem 1}, Lemma
\ref{lemma 1.1} and Propositions \ref{proposition 4},
\ref{proposition 5}, the minimal number of generators for
$Lift(\varphi_k)$
can be calculated as follows:
\begin{eqnarray*}
{ } & { } & (2k-1)\cdot \left(
\begin{array}{c}
2k-1 \\
1
\end{array}
\right)
-
\left({((2k-1)-(2k-2))}\cdot {}_{1}\delta(\varphi_k)+{}_{1}\gamma(\varphi_k)-{}_0\gamma(\varphi_k)\right)  \\
{ } & = &
(2k-1)^2-\left(((2k-1)-(2k-2))\cdot (2k-2)\cdot k\right. \\
{ } & { } & \left. \qquad\qquad\qquad\qquad +(2k-2)\cdot (k-1)-(k-1)\right)  \\
{ } & = & 3k-2.
\end{eqnarray*}
It has been verified in \cite{hollandmond} that the minimal number
of generators for $Lift(\varphi_k)$
is exactly $3k-2$ in the complex case and in the case a set of
generators has been obtained in \cite{houstonlittlestone} (see also
\cite{brucewest}). In Subsection \ref{subsection 6.1}, a set of
linear parts of generators for $Lift(\varphi_k)$ is obtained by our
method for any $k\ge 2$. }
\end{example}
\begin{example}\label{example 3.3}
{\rm Let $\psi_n: (\mathbb{K}^n,0)\to (\mathbb{K}^{2n-1},0)$ be
given by
\[
\psi_n(v_1, \ldots, v_{n-1},y)=(v_1, \ldots, v_{n-1}, y^2,  v_1y,
\ldots, v_{n-1}y).
\]
Then, it is known that $f$ is an isolated stable mono-germ by
\cite{whitney1} or \cite{whitney2} or \cite{mather5}. Thus, by
Proposition \ref{proposition new 1}, ${}_0\overline{\omega}\psi_n$
is bijective . Therefore, it follows that
$i_1(\psi_n)=i_2(\psi_n)=0$. By Theorem \ref{theorem 1}, Lemma
\ref{lemma 1.1} and Propositions \ref{proposition 4},
\ref{proposition 5}, the minimal number of generators for
$Lift(\psi_n)$
can be calculated as follows:
\begin{eqnarray*}
{ } & { } & (2n-1)\cdot \left(
\begin{array}{c}
2n-1 \\
1
\end{array}
\right)
-
\left(((2n-1)-n)\cdot {}_{1}\delta(\varphi)+{}_{1}\gamma(\varphi)-{}_0\gamma(\varphi)\right)  \\
{ } & = &
(2n-1)^2-\left((n-1)\cdot n\cdot 2+n\cdot (2-1)-(2-1)\right)  \\
{ } & = & 2n^2-3n+2.
\end{eqnarray*}
In the case that $n=2$, $\psi_2$ equals $\varphi_2$ of Example
\ref{example 3.2}. Thus, in this case, It has been verified in
\cite{brucewest} and \cite{hollandmond} that the minimal number of
generators for $Lift(\psi_n)$
is exactly $4$ in the complex case and a set of generators has been
obtained in \cite{brucewest} and \cite{houstonlittlestone}. In
Subsection \ref{subsection 6.2}, a set of generators for
$Lift(\psi_n)$ is obtained by our method for any $n\ge 2$. }
\end{example}
\begin{example}\label{example 3.4}
{\rm Examples \ref{example 3.1}, \ref{example 3.2} and \ref{example
3.3} can be generalized as follows.    Let $f: (\mathbb{K},0)\to
(\mathbb{K}^p,0)$ $(p\ge 2)$ be an \textcolor{black}{analytic}
map-germ such that $2\le \delta(f)<\infty$ and let $F:
(\mathbb{K}\times \mathbb{K}^c,0)\to
(\mathbb{K}^p\times\mathbb{K}^c,0)$ be a $\mathcal{K}$-miniversal
unfolding of $f$, where $\mathcal{K}$-{\it miniversal unfolding of}
$f$ is a map-germ given by (5.8) of \cite{mather4} with $c=r$. Then,
by \cite{mather4} or \cite{mather5}, $F$ is an isolated stable
mono-germ. Thus, by Proposition \ref{proposition new 1},
${}_0\overline{\omega}F$ is bijective. Note that $c=p\delta(f)-1-p$
by theorem 4.5.1 of \cite{wall}. By Theorem \ref{theorem 1}, Lemma
\ref{lemma 1.1} and Propositions \ref{proposition 4},
\ref{proposition 5}, the minimal number of generators for $Lift(F)$
can be calculated as follows:
\begin{eqnarray*}
{ } & { } & (p+c)\cdot \left(
\begin{array}{c}
p+c \\
1
\end{array}
\right)
-
\left(\left((p+c)-(1+c)\right)\cdot {}_{1}\delta(F)+{}_{1}\gamma(F)-{}_0\gamma(F)\right)  \\
{ } & = &
(p+c)^2- \left((p-1)\cdot (1+c)\cdot \delta(f)+c\cdot \left(\delta(f)-1\right) \right) \\
{ } & = & p^2\cdot \delta(f)-p\cdot \delta(f)+\delta(f)-p.
\end{eqnarray*}
By Mather's classification theorem (theorem A of \cite{mather4}),
proposition (1.6) of \cite{mather4}, Mather's normal form theorem
for a stable map-germ (theorem (5.10) of \cite{mather4}), the fact
that the sharp inequality $p^2\delta(f)-p\delta(f)+\delta(f)-p>p+c$
holds (since $p, \delta(f)\ge 2$)
and the fact that the module of liftable vector fields over an
immersive stable multigerm is a free module if and only if $p=n+1$,
we have the following:
\begin{proposition}\label{proposition 6}
Let $f: (\mathbb{K}^n,S)\to (\mathbb{K}^p,0)$ $(n<p)$ be a stable
multigerm of corank at most one. Then, $Lift(f)$
is a free module if and only if the properties $p=n+1$ and
$\delta(f)=|S|$ are satisfied.
\end{proposition}
}
\end{example}
\begin{example}\label{example multistable}
{\rm Let $f: (\mathbb{K}^2,S)\to (\mathbb{K}^2,0)$ be given by
$(x,y)\mapsto (x,y^2)$, $(x,y)\mapsto (x^2, y)$. Then, it is known
that $f$ is an isolated multigerm by \cite{whitney3} or
\cite{mather5}. Thus, by Proposition \ref{proposition new 1},
${}_0\overline{\omega}f$ is bijective. Therefore, it follows that
$i_1(f)=i_2(f)=0$. By Theorem \ref{theorem 1}, Lemma \ref{lemma 1.1}
and Propositions \ref{proposition 4}, \ref{proposition 5}, the
minimal number of generators for $Lift(f)$
is the following:
\begin{eqnarray*}
{ } & { } &
\dim_{\mathbb{K}}\frac{m_0\theta_0(2)}{m_0^{2}\theta_0(2)}-
\left((2-2){}_{1}\delta(f)+{}_{1}\gamma(f)-{}_0\gamma(f)\right)  \\
{ } & = &
2^2-\left((2-2)\times 2\times 4+2\times (4-2)-(4-2)\right)  \\
{ } & = & 2.
\end{eqnarray*}
In this case we can construct easily a basis of $Lift(f)$ consisting
of $2$ vector fields (see Subsection \ref{subsection 6.4}). }
\end{example}
\begin{example}\label{example 3.5}
{\rm Let $f: (\mathbb{K},S)\to (\mathbb{K}^2,0)$ be any one of the
following three.
\begin{enumerate}
\item $x\mapsto (x^4, x^5+x^7)$ (taken from \cite{brucegaffney}).
\item $x\mapsto (x^2,x^3), x\mapsto (x^3, x^2)$ (taken from \cite{kolgushkinsadykov}).
\item $x\mapsto (x,0), x\mapsto (0,x), x\mapsto (x^2,x^3+x^4)$ (taken from \cite{kolgushkinsadykov}).
\end{enumerate}
It has been shown in \cite{brucegaffney} or \cite{kolgushkinsadykov}
that $T\mathcal{K}(f)=T\mathcal{A}(f)$ is satisfied. Thus, by
Proposition \ref{proposition new 1}, ${}_1\overline{\omega}f$ is
surjective. We can confirm easily that the following equality holds.
\[
2\cdot \left(
\begin{array}{c}
2 \\
1
\end{array}
\right)
= (2-1)\cdot {}_{1}\delta(f)+{}_{1}\gamma(f)-{}_0\gamma(f).
\]
Thus, ${}_1\overline{\omega}f$ is injective by Proposition
\ref{proposition 4}. Therefore, it follows that $i_1(f)=i_2(f)=1$.
By Theorem \ref{theorem 1}, Lemma \ref{lemma 1.1} and Propositions
\ref{proposition 4}, \ref{proposition 5}, the minimal number of
generators for $Lift(f)$
can be calculated as follows:
\begin{eqnarray*}
{ } & { } & 2\cdot \left(
\begin{array}{c}
3 \\
2
\end{array}
\right)
-
\left({(2-1)}\cdot {}_{2}\delta(f)+{}_{2}\gamma(f)-{}_1\gamma(f)\right)  \\
{ } & = &
2\cdot 3-\left((2-1)\cdot 1\cdot 4+(4-|S|)-(4-|S|)\right)  \\
{ } & = & 2.
\end{eqnarray*}
In the case $\mathbb{K}=\mathbb{C}$, it has been known that any
plane algebraic curve is a free divisor by \cite{saito}. Thus, by
combining \cite{damon} and \cite{saito}, it has been known that the
minimal number of generators for $Lift(f)$
is $2$ in the complex case. }
\end{example}
\begin{example}\label{example 3.6}
{\rm Let $f: (\mathbb{K}^2,0)\to (\mathbb{K}^2,0)$ be given by
$f(x,y)=(x,xy+y^5\pm y^7)$ (taken from \cite{rieger}). It has been
shown in \cite{rieger} that $T\mathcal{K}(f)=T\mathcal{A}(f)$ is
satisfied. Thus, by Proposition \ref{proposition new 1},
${}_1\overline{\omega}f$ is surjective. It is easily seen that the
following equality holds.
\[
2\cdot \left(
\begin{array}{c}
2 \\
1
\end{array}
\right)
= (2-2)\cdot {}_{1}\delta(f)+{}_{1}\gamma(f)-{}_0\gamma(f).
\]
Thus, ${}_1\overline{\omega}f$ is injective by Proposition
\ref{proposition 4}. Therefore, it follows that $i_1(f)=i_2(f)=1$.
By Theorem \ref{theorem 1}, Lemma \ref{lemma 1.1} and Propositions
\ref{proposition 4}, \ref{proposition 5}, the minimal number of
generators for $Lift(f)$
can be calculated as follows.
\begin{eqnarray*}
{ } & { } & 2\cdot \left(
\begin{array}{c}
3 \\
2
\end{array}
\right)
-
\left((2-2)\cdot {}_{2}\delta(f)+{}_{2}\gamma(f)-{}_1\gamma(f)\right)  \\
{ } & = &
2\cdot 3-\left((2-2)\cdot 3\cdot 5+3\cdot (5-1)-2\cdot (5-1)\right)  \\
{ } & = & 2.
\end{eqnarray*}
As same as Example \ref{example 3.5}, it has been known that the
minimal number of generators for $Lift(f)$
is $2$ in the complex case. }
\end{example}
\begin{example} \label{example riegerruaswikatique}
{\rm Let $f: (\mathbb{K}^4,0)\to (\mathbb{K}^5,0)$ be given by the
following:
\[
f(x_1, x_2, x_3, y)=(x_1, x_2, x_3, y^4+x_1y, y^6+y^7+x_2y+x_3y^2).
\]
This example is taken from \cite{riegerruaswikatique} where the
property $T\mathcal{K}(f)=T\mathcal{A}(f)$ has been shown. Thus, by
Proposition \ref{proposition new 1}, ${}_1\overline{\omega}f$ is
surjective. It is easily seen that the following equality holds.
\[
5\cdot \left(
\begin{array}{c}
5 \\
1
\end{array}
\right)
= (5-4)\cdot {}_{1}\delta(f)+{}_{1}\gamma(f)-{}_0\gamma(f).
\]
Thus, ${}_1\overline{\omega}f$ is injective by Proposition
\ref{proposition 4}. Therefore, it follows that $i_1(f)=i_2(f)=1$.
By Theorem \ref{theorem 1}, Lemma \ref{lemma 1.1} and Propositions
\ref{proposition 4}, \ref{proposition 5}, the minimal number of
generators for $Lift(f)$
can be calculated as follows:
\begin{eqnarray*}
{ } & { } & 5\cdot \left(
\begin{array}{c}
6 \\
2
\end{array}
\right)
-
\left((5-4)\cdot {}_{2}\delta(f)+{}_{2}\gamma(f)-{}_1\gamma(f)\right)  \\
{ } & = &
5\cdot 15-\left((5-4)\cdot 10\cdot 4+10\cdot (4-1)-4\cdot (4-1)\right)  \\
{ } & = & 17.
\end{eqnarray*}
}
\end{example}
\begin{example}\label{example multicusp}
{\rm Let $c: \mathbb{K}\to \mathbb{K}^2$ be the map defined by
$c(x)=(x^2, x^3)$. For any real number $\theta$, we let $R_{\theta}:
\mathbb{K}^2\to \mathbb{K}^2$ be the rotation of $\mathbb{K}^2$
about the origin with respect to the angle $\theta$:
$$
R_{\theta} \left(
\begin{array}{c}
X \\
Y
\end{array}
\right) = \left(
\begin{array}{rr}
\cos\theta & -\sin\theta \\
\sin\theta & \cos\theta
\end{array}
\right) \left(
\begin{array}{c}
X \\
Y
\end{array}
\right).
$$
For any $i\in \mathbb{N}$, let $\theta_0, \ldots, \theta_i$ be real
numbers such that $0\le \theta_j<2\pi$ $(0\le j\le i)$ and $0\ne
|\theta_j-\theta_k|\ne \pi$ $(j\ne k)$. Set $S=\{s_0, \ldots, s_i\}$
$(s_j\ne s_k \mbox{ if }j\ne k)$ and define $c_{\theta_j}:
(\mathbb{K}, s_j)\to (\mathbb{K}^2,0)$ as
$c_{\theta_j}(x)=R_{\theta_j}\circ c(x_j)$, where $x_j=x-s_j$.
A multigerm $\{c_{\theta_0}, \ldots, c_{\theta_i}\}:
(\mathbb{K},S)\to (\mathbb{K}^2,0)$, which is called a {\it
multicusp},
is denoted by $c_{(\theta_0, \ldots, \theta_i)}$.
\par
In \cite{mizotanishimura}, it has been shown that $i_1(c_{(\theta_0,
\ldots, \theta_i)})=i_2(c_{(\theta_0, \ldots, \theta_i)})=i$ for any
$i\in \mathbb{N}$. Thus,
by Theorem \ref{theorem 1}, Lemma \ref{lemma 1.1} and Propositions
\ref{proposition 4}, \ref{proposition 5}, the minimal number of
generators for $Lift(c_{(\theta_0, \ldots, \theta_i)})$
can be calculated as follows.
\begin{eqnarray*}
{ } & { } & 2\cdot \left(
\begin{array}{c}
2+i \\
i+1
\end{array}
\right)
- \left((2-1)\cdot {}_{i+1}\delta(c_{(\theta_0, \ldots, \theta_i)})+
{}_{i+1}\gamma(c_{(\theta_0, \ldots, \theta_i)})-{}_i\gamma(c_{(\theta_0, \ldots, \theta_i)})\right)  \\
{ } & = & 2\cdot (2+i)-\left((2-1)\cdot 1\cdot \delta(c_{(\theta_0,
\ldots, \theta_i)}) +1\cdot \gamma(c_{(\theta_0, \ldots,
\theta_i)})-1\cdot
\gamma(c_{(\theta_0, \ldots, \theta_i)})\right)  \\
{ } & = & 2.
\end{eqnarray*}
As same as Example \ref{example 3.5}, it has been already known that
the minimal number of generators for $Lift(c_{(\theta_0, \ldots,
\theta_i)})$
is $2$ in the complex case. }
\end{example}
\section{Proof of Theorem \ref{theorem 1}}
\label{section 4}
Since ${}_i\overline{\omega}f$ is surjective, by Lemma 1.1 we have
that ${}_j\overline{\omega}f$ is surjective for any $j>i$. Since
${}_i\overline{\omega}f$ is injective, any $\eta\in \theta_0(p)$
such that $\omega f(\eta)\in T\mathcal{R}_e(f)$ is contained in
$m_0^{i+1}\theta_0(p)$. Set $\rho(f)= \dim_{\mathbb{K}}\mbox{\rm
ker}({}_{i+1}\overline{\omega}f)$. Then, since
${}_i\overline{\omega}f$ is bijective, $\rho(f)$ must be positive by
Corollary 1. Let $\{\eta_1+m_0^{i+2}\theta_0(p), \ldots,
\eta_{\rho(f)}+m_0^{i+2}\theta_0(p)\}$ be a basis of
$\mbox{\rm ker}({}_{i+1}\overline{\omega}f)$. Then, we have that
\[
\eta_j\circ f\in T\mathcal{R}_e(f)\cap
f^*m_0^{i+1}\theta_S(f)+f^*m_0^{i+2}\theta_S(f)\quad (1\le j\le
\rho(f)).
\]
Since ${}_{i+2}\overline{\omega}f$ is surjective, we have the
following:
\begin{eqnarray*}
{ } & { } & T\mathcal{R}_e(f)\cap f^*m_0^{i+1}\theta_S(f)+f^*m_0^{i+2}\theta_S(f) \\
{ } & = & T\mathcal{R}_e(f)\cap f^*m_0^{i+1}\theta_S(f)+
T\mathcal{R}_e(f)\cap f^*m_0^{i+2}\theta_S(f)+\omega
f(m_0^{i+2}\theta_0(p)).
\end{eqnarray*}
Thus, for any $j$ $(1\le j\le \rho(f))$ there exists
$\widetilde{\eta}_j\in m_0^{i+2}\theta_0(p)$ such that
$(\eta_j+\widetilde{\eta}_j)\circ f\in T\mathcal{R}_e(f)\cap
T\mathcal{L}_e(f)$. Let $A$ be the $C_0$-module generated by
$\eta_j+\widetilde{\eta}_j$ $(1\le j\le \rho(f))$.
\par
Let $\hat{\omega}f: \theta_0(p)\to
\frac{\theta_S(f)}{T\mathcal{R}_e(f)}$ be given by
$\hat{\omega}f(\eta)=\omega f(\eta)+T\mathcal{R}_e(f)$. Then,
$\mbox{\rm ker}(\hat{\omega}f)$ is the set of vector fields liftable
over $f$.
In order to show that $\mbox{\rm ker}(\hat{\omega}f)=A$, we consider
the following commutative diagram.
\[
\begin{CD}
{ } @. 0 @>>> 0 @>>>\mbox{\rm ker}(b_3) \\
@. @VVV @VVV @VVV \\
{ } @. m_0\mbox{\rm ker}(\hat{\omega}f)   @>{a_2}>>
m_0^{i+2}\theta_0(p) @>{c_2}>>
\frac{f^*m_0^{i+2}\theta_S(f)}{f^*m_0\left(T\mathcal{R}_e(f)\cap
f^*m_0^{i+1}\theta_S(f)\right)}
@>>> 0 \\
@. @VV{b_1}V @VV{b_2}V @VV{b_3}V \\
0 @>>> \mbox{\rm ker}(\hat{\omega}f)   @>{a_1}>>
m_0^{i+1}\theta_0(p) @>{c_1}>>
\frac{f^*m_0^{i+1}\theta_S(f)}{T\mathcal{R}_e(f)\cap
f^*m_0^{i+1}\theta_S(f)}
@.  { }  \\
@. @VVV @VVV @VVV \\
{ } @. \mbox{\rm coker}(b_1) @>{d_1}>>\mbox{\rm coker}(b_2)
@>{d_2}>>\mbox{\rm coker}(b_3)
\end{CD}
\]
Here, $a_j$ $(j=1, 2)$, $b_j$ $(j=1, 2)$ are inclusions, $b_3$ is
defined by $b_3([\eta]_{i+2})=[\eta]_{i+1}$ and $c_j$ $(j=1, 2)$ are
defined by $c_j(\eta)=[\omega f(\eta)]_{i+j}$, where
$[\eta]_{i+1}=\eta +T\mathcal{R}_e(f)\cap f^*m_0^{i+1}\theta_S(f)$
and $[\eta]_{i+2}=\eta +f^*m_0\left(T\mathcal{R}_e(f)\cap
f^*m_0^{i+1}\theta_S(f)\right)$.
\begin{lemma}
\[
f^*m_0\left(T\mathcal{R}_e(f)\cap
f^*m_0^{i+1}\theta_S(f)\right)=T\mathcal{R}_e(f)\cap
f^*m_0^{i+2}\theta_S(f).
\]
\end{lemma}
\noindent \underline{\it Proof of Lemma 4.1.}\quad It is clear that
$f^*m_0\left(T\mathcal{R}_e(f)\cap f^*m_0^{i+1}\theta_S(f)\right)
\subset T\mathcal{R}_e(f)\cap f^*m_0^{i+2}\theta_S(f)$. Thus, in the
following we concentrate on showing its converse. Let $\eta$ be an
element of $T\mathcal{R}_e(f)\cap f^*m_0^{i+2}\theta_S(f)$. Let
$f_j$ be a branch of $f$, namely, $f_j=f|_{(\mathbb{K}^n,s_j)}$
$(1\le j\le |S|)$. Then, we have that $\eta\in
T\mathcal{R}_e(f_j)\cap f_j^*m_0^{i+2}\theta_{s_j}(f_j)$ for any $j$
$(1\le j\le |S|)$. Thus, there exists $\xi_j\in \theta_{s_j}(n)$
such that $tf_j(\xi_j)=\eta$.
\par
Since $f$ is of corank at most one, for any $j$ $(1\le j\le |S|)$
there exist germs of diffeomorphism $h_j: (\mathbb{K}^n,s_j)\to
(\mathbb{K}^n,s_j)$ and  $H_j: (\mathbb{K}^p,0)\to (\mathbb{K}^p,0)$
such that $H_j\circ f_j\circ h_j^{-1}$ has the following form:
\begin{eqnarray*}
{ } & { } & H_j\circ f_j\circ h_j^{-1}(x, y) \\
{ } & = & (x, y^{\delta(f_j)}+f_{j, n}(x, y), f_{j, n+1}(x, y),
\ldots, f_{j, p}(x, y)).
\end{eqnarray*}
Here, $x$ stands for $(x_1, \ldots, x_{n-1})$ and $x_1, \ldots,
x_{n-1}, y$ are local coordinates of the coordinate system $(U_j,
h_j)$ at $s_j$ and $f_{j, q}$ satisfies $f_{j,q}(0, \ldots, 0,
y)=o(y^{\delta(f_j)})$ for any $q$ $(n\le q\le p)$. Set
\[
\xi_j=\sum_{m=1}^{n-1}\xi_{j, m}\frac{\partial}{\partial
x_m}+\xi_{j, n}\frac{\partial}{\partial y} \mbox{ and }
\eta=\sum_{q=1}^p \eta_q \frac{\partial}{\partial X_q}.
\]
Then, by the above form of $H_j\circ f_j\circ h_j^{-1}$ and the
equality $tf_j(\xi_j)=\eta$, the following hold:
\begin{eqnarray}
& { } & \xi_{j, m}(x_1, \ldots, x_{n-1},y)  =  \eta_m(x_1, \ldots, x_{n-1},y) \quad (1\le m\le n-1)\\
& { } & \lambda(x_1, \ldots, x_{n-1},y)\xi_{j, n}(x_1, \ldots,
x_{n-1},y)=\mu(x_1, \ldots, x_{n-1},y),
\end{eqnarray}
where $\lambda=\delta(f_j)y^{\delta(f_j)-1} +\frac{\partial
f_{j,n}}{\partial y}$ and
$\mu=\eta_n-\sum_{m=1}^{n-1}\xi_{j,m}\frac{\partial
f_{j,n}}{\partial x_m}$. Since $\eta_q\in f^*m_0^{i+2}C_{s_j}$ for
any $q$ $(1\le q\le p)$, by (4.1) we have that $\xi_{j, m}\in
f^*m_0^{i+2}C_{s_j}$ for any $m$ $(1\le m\le n-1)$.
\par
Since $f_{j,q}(0, \ldots, 0, y)=o(y^{\delta(f_j)})$ for any $q$
$(n\le q\le p)$, we have the following properties:
\begin{enumerate}
\item $Q(f_j)=Q(x_1, \ldots, x_{n-1}, y^{\delta(f_j)})=Q(x_1, \ldots, x_n, y^\lambda)$.
\item $[1], [y], \ldots, [y^{\delta(f_j)-2}], [\lambda]$ constitute a basis of $Q(f_j)$.
\end{enumerate}
Thus, by the preparation theorem, $C_{s_j}$ is generated by $1, y,
\ldots, y^{\delta(f_j)-2}, \lambda$ as $C_0$-module via $f_j$.
Therefore, for any positive integer $r$, $f_j^*m_0^{r}C_{s_j}$ is
generated by elements of the union of the following three sets $U_r,
V_r, W_r$ as $C_0$-module via $f_j$.
\begin{eqnarray*}
U_{r} & = & \left\{x_1^{k_1}\cdots
x_{n-1}^{k_{n-1}}\lambda^{k_n}y^{k_n+\ell}\; \left|\;
k_m\ge 0, \sum_{m=1}^{n-1}k_m=r-k_n<r, 0 \le \ell \le \delta(f_j)-2\right.\right\}, \\
V_{r} & = & \left\{x_1^{k_1}\cdots
x_{n-1}^{k_{n-1}}\lambda^{k_n+1}y^{k_n}\; \left|\;
k_m\ge 0, \sum_{m=1}^{n-1}k_m=r-k_n, \right.\right\}, \\
W_{r} & = & \left\{x_1^{k_1}\cdots x_{n-1}^{k_{n-1}}y^{\ell}\;
\left|\; k_m\ge 0, \sum_{m=1}^{n-2}k_m=r, 0\le \ell \le
\delta(f_j)-2\right.\right\}.
\end{eqnarray*}
Then, by using these notations,
for any $m$ $(1\le m\le n-1)$ $\xi_{j, m}$ can be expressed as
follows:
 \[
 \xi_{j, m}=
\sum_{u\in U_{i+2}}\varphi_{u, j,m}u + \sum_{v\in
V_{i+2}}\varphi_{v, j,m}v + \sum_{w\in W_{i+2}}\varphi_{w, j,m}w,
\]
where $\varphi_{u, j, m}, \varphi_{v, j, m}, \varphi_{w, j, m}$ are
some elements of $C_{s_j}$.
\par
Next, we investigate $\xi_{j,n}$. Since $\mu$ has the form
$\mu=\eta_n-\sum_{m=1}^{n-1}\xi_{j,m}\frac{\partial
f_{j,n}}{\partial x_m}$ and $\eta_n, \xi_{j,m}$ are contained in
$f^*m_0^{i+2}C_{s_j}$, $\mu$ is contained in $f^*m_0^{i+2}C_{s_j}$.
On the other hand, $\lambda$ must divide $\mu$ by (4.2). Thus, $\mu$
is generated by elements of $U_{i+2}\cup V_{i+2}$. Hence,
$\xi_{j,n}=\frac{\mu}{\lambda}$ can be expressed as follows:
 \[
 \xi_{j, n}=
\sum_{u\in U_{i+2}}\varphi_{u, j,n}\frac{u}{\lambda} + \sum_{v\in
V_{i+2}}\varphi_{v, j,n}\frac{v}{\lambda},
\]
where $\varphi_{u, j, n}, \varphi_{v, j, n}$ are some elements of
$C_{s_j}$. Since $\frac{u}{\lambda}\in U_{i+1}\cup V_{i+1}\cup
W_{i+1}$ for any $u\in U_{i+2}$ and $\frac{v}{\lambda}\in U_{i+1}$
(resp., $\frac{v}{\lambda}\in V_{i+1}$) if $\delta(f_j)\ge 2$
(resp., $\delta(f_j)=1$) for any $v\in V_{i+2}$, $\xi_{j,n}$ is
belonging to $f_j^*m_0^{i+1}C_{s_j}$.
\par
\smallskip
Since $f_j^*m_0^{i+2}C_{s_j}\subset f_j^*m_0^{i+1}C_{s_j}$, for any
$j$ $(1\le j\le |S|)$ and any $m$ $(1\le m\le n-1)$, we have the
following:
\begin{eqnarray*}
{ } & { } & tf_j\left(\xi_{j,m}\frac{\partial}{\partial x_m}\right) \\
{ } & = & tf_j\left(\left(\sum_{u\in U_{i+1}}
\widetilde{\varphi}_{u,j,m}u+ \sum_{v\in V_{i+1}}
\widetilde{\varphi}_{v,j,m}v +
\sum_{w\in W_{i+1}} \widetilde{\varphi}_{w,j,m}w\right)\frac{\partial}{\partial x_m} \right) \\
{ } & = & \sum_{u\in U_{i+1}} u\left(
tf_j\left(\widetilde{\varphi}_{u,j,m}\frac{\partial}{\partial
x_m}\right)\right) +
\sum_{v\in V_{i+1}} v\left( tf_j\left(\widetilde{\varphi}_{v,j,m}\frac{\partial}{\partial x_m}\right)\right) \\
{ } & { } & \qquad\qquad + \sum_{w\in W_{i+1}} w\left(
tf_j\left(\widetilde{\varphi}_{w,j,m}\frac{\partial}{\partial
x_m}\right)\right),
\end{eqnarray*}
where $\widetilde{\varphi}_{u, j, m}, \widetilde{\varphi}_{v, j, m},
\widetilde{\varphi}_{w, j, m}$ are some elements of $C_{s_j}$.
Moreover, for any $j$ $(1\le j\le |S|)$ we have the following:
\begin{eqnarray*}
{ } & { } & tf_j\left(\xi_{j,n}\frac{\partial}{\partial y}\right) \\
{ } & = & tf_j\left(\left(\sum_{u\in U_{i+1}} \psi_{u,j,n}u+
\sum_{v\in V_{i+1}} \psi_{v,j,n}v+
\sum_{w\in W_{i+1}} \psi_{w,j,n}w\right)\frac{\partial}{\partial y}\right) \\
{ } & = & \sum_{u\in U_{i+1}} u
\left(tf_j\left(\psi_{u,j,n}\frac{\partial}{\partial
y}\right)\right) +
\sum_{v\in V_{i+1}} v \left(tf_j\left(\psi_{v,j,n}\frac{\partial}{\partial y}\right)\right) \\
{ } & { } & \qquad\qquad + \sum_{w\in W_{i+1}} w
\left(tf_j\left(\psi_{w,j,n}\frac{\partial}{\partial
y}\right)\right),
\end{eqnarray*}
where $\psi_{u, j, n}, \psi_{v, j, n}, \psi_{w, j, n}$ are elements
of $C_{s_j}$. Since the union $U_{i+1}\cup V_{i+1}\cup W_{i+1}$ is a
finite set,
we have the following:
\[
\eta=tf_j\left(\sum_{m=1}^{n-1}\xi_{j,m}\frac{\partial}{\partial
x_m}+\xi_{j,n}\frac{\partial}{\partial y}\right) \in
f_j^*m_0^{i+1}\left(T\mathcal{R}_e(f_j)\cap
f_j^*m_0\theta_{s_j}(f_j)\right).
\]
Notice that $i+1\ge 1$ and $f_j$ is any branch of $f$.      Hence,
we have that
$$
\eta\in f^*m_0\left(T\mathcal{R}_e(f)\cap
f^*m_0^{i+1}\theta_S(f)\right).
$$
  \hfill $\Box$
\par
\smallskip
\noindent Lemma 4.1 implies that $c_2$ is surjective, thus even the
second row sequence is exact. Lemma 4.1 implies also that $b_3$ is
injective
and thus $\mbox{\rm ker}(b_3)=0$. Hence, by the snake lemma, we see
that $d_1$ is injective. On the other hand, since
there exists an isomorphism
\[
\varphi: \frac{f^*m_0^{i+1}\theta_S(f)}{T\mathcal{R}_e(f)\cap
f^*m_0^{i+1}\theta_S(f)+f^*m_0^{i+2}\theta_S(f)} \to \mbox{\rm
coker}(b_3)
\]
such that $d_2=\varphi\circ  {}_{i+1}\overline{\omega}f$ we have
that $\mbox{\rm ker}(d_2) = \mbox{\rm ker}(\varphi\circ
{}_{i+1}\overline{\omega}f) = \mbox{\rm
ker}({}_{i+1}\overline{\omega}f)$. Therefore, we have the following:
\[
\dim_{\mathbb{K}}\frac{\mbox{\rm ker}(\hat{\omega}f)}{m_0\mbox{\rm
ker}(\hat{\omega}f)} = \dim_{\mathbb{K}}\mbox{\rm
ker}({}_{i+1}\overline{\omega}f) =\rho(f)
=\dim_{\mathbb{K}}\frac{A}{m_0A}.
\]
Moreover, $A$ is a submodule of $\mbox{\rm ker}(\hat{\omega}f)$
\textcolor{black}{and our category is the
\textcolor{black}{analytic} category}. Therefore, we have that
$\mbox{\rm ker}(\hat{\omega}f)=A$. \hfill Q.E.D.
\section{Proof of Theorem \ref{theorem 2}}\label{section 5}
We first show the assertion (1) of Theorem \ref{theorem 2}. Let
$F(x, \lambda)=(f_\lambda(x), \lambda)$ $(f_0=f)$ be a stable
unfolding of $f$, Since $\eta$ is an element of $Lift(F)$, by
definition, there exists a vector field $\xi\in \theta_S(n+r)$ such
that $dF\circ \xi=\eta\circ F$. Set $x=(x_1, \ldots, x_n)$,
$\lambda=(\lambda_1, \ldots, \lambda_r)$, $f_\lambda=(f_{\lambda,
1}, \ldots, f_{\lambda, p})$ and $\xi(x, \lambda)=(\xi_1(x,\lambda),
\ldots, \xi_{n+r}(x, \lambda))$. Then, we have the following:
\[
\left(
\begin{array}{cc}
\frac{\partial f_{\lambda}}{\partial x}(x, \lambda) & \frac{\partial f_{\lambda}}{\lambda}(x, \lambda) \\
0 & E_r
\end{array}
\right) \left(
\begin{array}{c}
\xi_1(x, \lambda) \\
\vdots \\
\xi_n(x, \lambda) \\
\xi_{n+1}(x, \lambda) \\
\vdots \\
\xi_{n+r}(x, \lambda)
\end{array}
\right) = \left(
\begin{array}{c}
\eta_1(F(x, \lambda)) \\
\vdots \\
\eta_p(F(x, \lambda)) \\
\eta_{p+1}(F(x, \lambda)) \\
\vdots \\
\eta_{p+r}(F(x, \lambda))
\end{array}
\right).
\]
Here, $\frac{\partial f_{\lambda}}{\partial x}(x, \lambda)$ is the
$p\times n$ matrix whose $(i,j)$ elements is $\frac{\partial
f_{\lambda, i}}{\partial x_j}(x, \lambda)$, $\frac{\partial
f_{\lambda}}{\partial \lambda}(x, \lambda)$ is the $p\times r$
matrix whose $(i,k)$ elements is $\frac{\partial f_{\lambda,
i}}{\partial \lambda_k}(x, \lambda)$, and $E_r$ stands for the
$r\times r$ unit matrix. In particular, we have the following two:
\[
\left( \frac{\partial f_{\lambda}}{\partial x}(x, 0) \right) \left(
\begin{array}{c}
\xi_1(x, 0) \\
\vdots \\
\xi_n(x, 0)
\end{array}
\right) + \left( \frac{\partial f_{\lambda}}{\partial \lambda}(x, 0)
\right) \left(
\begin{array}{c}
\xi_{n+1}(x, 0) \\
\vdots \\
\xi_{n+r}(x, 0)
\end{array}
\right) = \left(
\begin{array}{c}
\eta_1(F(x, 0)) \\
\vdots \\
\eta_p(F(x, 0))
\end{array}
\right)
\]
and
\[
\xi_{n+k}(x,0)=\eta_{p+k}(F(x,0))=\eta_{p+k}(f(x),0)=0\quad (1\le
k\le r). \leqno{(*)}
\]
The last equality of $(*)$ is obtained from the assumption that
$\eta\in Lift(g)$ and the fact that $Lift(g)$ is generated by
$\frac{\partial}{\partial X_1}, \ldots, \frac{\partial}{\partial
X_p}$
 and $\Lambda_1\frac{\partial}{\partial \Lambda_1}, \ldots, \Lambda_r\frac{\partial}{\partial \Lambda_r}$.
Thus, we have the following:
\[
\left( \frac{\partial f_{\lambda}}{\partial x}(x, 0) \right) \left(
\begin{array}{c}
\xi_1(x, 0) \\
\vdots \\
\xi_n(x, 0)
\end{array}
\right) = \left(
\begin{array}{c}
\eta_1(F(x, 0)) \\
\vdots \\
\eta_p(F(x, 0))
\end{array}
\right).
\]
Therefore, $\overline{\eta}$ is a liftable vector field of $f$.
\par
\medskip
We next show the assertion (2) of Theorem \ref{theorem 2}. Since
$\overline{\eta}\in Lift(f)$, by definition, there exists a vector
field $\overline{\xi}\in \theta_S(n)$ such that
$df(\overline{\xi})=\overline{\eta}\circ f$. Set $\eta(X,
\Lambda)=(\overline{\eta}(X), 0)$ and $\xi(x,
\lambda)=(\overline{\xi}(x), 0)$.    Set also
$\widetilde{\eta}=\eta\circ F-dF(\xi)\in \theta_{S\times\{0\}}(F)$.
It is not difficult to see that $\widetilde{\eta}(x,0)=(0,0)$. Thus,
by the preparation theorem, there must exist a vector field
$\widetilde{\eta}_1\in \theta_{S\times \{0\}}(F)$ such that
$\widetilde{\eta}=\Lambda\widetilde{\eta}_1$.
\par
Since $F$ is stable,  there exist $\hat{\xi}\in \theta_{S\times
\{0\}}(n+1)$ and $\hat{\eta}\in \theta_{(0,0)}(p+1)$ such that
$\widetilde{\eta}_1=dF(\hat{\xi})+\hat{\eta}\circ F$. We therefore
have that
\[
\eta\circ F - dF(\xi)=\widetilde{\eta}=\Lambda\widetilde{\eta}_1 =
\lambda(dF(\hat{\xi})+\hat{\eta}\circ F)=
dF(\lambda\hat{\xi})+(\Lambda\hat{\eta})\circ F.
\]
It is clear that $\eta-\Lambda\hat{\eta}\in Lift(F)$. Moreover, we
have that the $(p+1)$ component of  $\eta-\Lambda\hat{\eta}$ is
$0-(\Lambda\hat{\eta})_{p+1}=-\Lambda(\hat{\eta})_{p+1}$ , which
implies that $\eta-\Lambda\widetilde{\eta}\in Lift(g_1)$. \hfill
Q.E.D.

\section{How to construct liftable vector fields
}\label{section 6}
In principle, the proof of Theorem \ref{theorem 1} provides how to
construct generators for the module of liftable vector fields over a
given finitely determined multigerm $f$ satisfying the assumption of
Theorem \ref{theorem 1}. In Subsections \ref{subsection
6.1}--\ref{subsection 6.5}, we examine it by several examples. In
Subsections \ref{subsection 6.6} and \ref{subsection 6.7}, as an
application of Thereom \ref{theorem 2}, we explain how to construct
all liftable vector fields over a given \textcolor{black}{analytic}
multigerm admitting a one-parameter stable unfolding by several
examples.
\subsection{$Lift(\varphi_k)$ for $\varphi_k(u_1, \ldots, u_{k-2}, v_1, \ldots, v_{k-1}, y)= \\(u_1, \ldots, u_{k-2}, v_1, \ldots, v_{k-1},
y^k+\sum_{i=1}^{k-2}u_iy^i, \sum_{i=1}^{k-1}v_iy^i)$}
\label{subsection 6.1}\qquad \\
Since the purpose of this subsection is to examine that the proof of
Theorem \ref{theorem 1} in principle provides how to construct
generators for the module of liftable vector fields,
in order to avoid just long calculations, in this subsection we
restrict ourselves to obtain only the linear parts of generators for
the module of liftable vector fields over $\varphi_k$ of Example
\ref{example 3.2}. Note that the linear parts of generators
themselves are already useful to obtain the best lower bound for the
$\mathcal{A}_e$ codimensions of some multigerms (see
\cite{raulcidinharoberta}).
\begin{definition}
{\rm Let $\eta_1, \eta_2$ be vector fields along $\varphi_k$
(namely, $\eta_1, \eta_2\in \theta_S(\varphi_k)$).
\begin{enumerate}
\item \textcolor{black}{We denote} $\eta_1\equiv \eta_2$ (mod $T\mathcal{R}_e(\varphi_k)$)) if
$\eta_1-\eta_2\in T\mathcal{R}_e(\varphi_k)$.
\item \textcolor{black}{We denote} $\eta_1\equiv \eta_2$
(mod
$T\mathcal{R}_e(\varphi_k)+\varphi_k^*m_0^2\theta_S(\varphi_k)$)) if
$\eta_1-\eta_2\in
T\mathcal{R}_e(\varphi_k)+\varphi_k^*m_0^2\theta_S(\varphi_k)$.
\end{enumerate}
}
\end{definition}
Let $(U_1, \ldots, U_{k-2}, V_1, \ldots, V_{k-1}, W_1, W_2)$ be the
standard coordinates of $\mathbb{K}^{2k-1}$. Along the proof of
Theorem \ref{theorem 1}, we first look for clues of linear terms of
liftable vector fields in
$\frac{\varphi_k^*m_0\theta_S(\varphi_k)}{\varphi_k^*m_0^2\theta_S(\varphi_k)}$.
From the form of the Jacobian matrix of $\varphi_k$ and since the
minimal number of generators is $(3k-2)$ by Example \ref{example
3.2}, we can guess that clues of linear terms of liftable vector
fields are the following $(3k-2)$ vector fields along $\varphi_k$:
\begin{eqnarray*}
{ } & { } & (W_1\circ \varphi_k)\frac{\partial}{\partial W_1}, \;
(W_2\circ \varphi_k)\frac{\partial}{\partial W_1}, \;
(W_2\circ \varphi_k)\frac{\partial}{\partial W_2}, \\
{ } & { } & y^i(W_1\circ \varphi_k)\frac{\partial}{\partial W_1}\quad (1\le i\le k-2), \\
{ } & { } & y^i(W_2\circ \varphi_k)\frac{\partial}{\partial W_1}\quad (1\le i\le k-2), \\
{ } & { } & y^i(W_2\circ \varphi_k)\frac{\partial}{\partial
W_2}\quad (1\le i\le k-1).
\end{eqnarray*}
\par
First we try to find a vector field $\eta\in m_0\theta_0(2k-1)$ such
that $W_1\frac{\partial}{\partial W_1}\ne \eta$ and $(W_1\circ
\varphi_k)\frac{\partial}{\partial W_1}\equiv \eta\circ \varphi_k$
(mod $T\mathcal {R}_e(\varphi_k))$ because
$W_1\frac{\partial}{\partial W_1}- \eta$ must be a liftable vector
field for such a $\eta$.
\begin{eqnarray*}
{ } & { } & (W_1\circ \varphi_k)\frac{\partial}{\partial W_1} \\
{ } & = & \left(y^k+\sum_{i=1}^{k-2}u_iy^i\right)\frac{\partial}{\partial W_1} \\
{ } & \equiv &
\left(-\frac{1}{k}\sum_{i=1}^{k-2}iu_iy^i+\sum_{i=1}^{k-2}u_iy^i\right)\frac{\partial}{\partial
W_1}
-\sum_{i=1}^{k-1}iv_iy^{i-1}\frac{\partial}{\partial W_2}\quad (\mbox{mod}\; T\mathcal{R}_e(\varphi_k)) \\
{ } & \equiv & -\sum_{i=1}^{k-2}\frac{k-i}{k}(U_i\circ
\varphi_k)\frac{\partial}{\partial U_i}
+ \sum_{i=2}^{k-1}i(V_i\circ \varphi_k)\frac{\partial}{\partial V_{i-1}}  \\
{ } & { } & \qquad \qquad \qquad \qquad +(V_1\circ
\varphi_k)\frac{\partial}{\partial W_2}\quad (\mbox{mod}\;
T\mathcal{R}_e(\varphi_k)).
\end{eqnarray*}
Thus, it follows that the following is a liftable vector fields over
$\varphi_k$ where $\widetilde{\eta}_1=0$.
\[
\eta_1+\widetilde{\eta}_1=\sum_{i=1}^{k-2}\frac{k-i}{k}U_i\frac{\partial}{\partial
U_i} - \sum_{i=2}^{k-1}iV_i\frac{\partial}{\partial V_{i-1}}
+W_1\frac{\partial}{\partial W_1} -V_1\frac{\partial}{\partial W_2}.
\]
\par
Secondly, we try to find a vector field $\eta\in m_0\theta_0(2k-1)$
such that $W_2\frac{\partial}{\partial W_1}\ne \eta$ and $(W_2\circ
\varphi_k)\frac{\partial}{\partial W_1}\equiv \eta\circ \varphi_k$
(mod $T\mathcal
{R}_e(\varphi_k))+\varphi_k^2m_0^2\theta(\varphi_k)$).
\begin{eqnarray*}
(W_2\circ \varphi_k)\frac{\partial}{\partial W_1} & = &
\left(\sum_{i=1}^{k-1}v_iy^i\right)\frac{\partial}{\partial W_1} \\
{ } & \equiv & -\sum_{i=1}^{k-2}(V_i\circ
\varphi_k)\frac{\partial}{\partial U_i}\quad (\mbox{mod}\;
T\mathcal{R}_e(\varphi_k) +\varphi_k^2m_0^2\theta(\varphi_k)).
\end{eqnarray*}
Thus, it follows that there exists a liftable vector field over
$\varphi_k$ having the following form, where $\widetilde{\eta}_2\in
m_0^2\theta_0(p)$.
\[
\eta_2+\widetilde{\eta}_2=\sum_{i=1}^{k-2}V_i\frac{\partial}{\partial
U_i} +W_2\frac{\partial}{\partial W_1} + \mbox{ higher terms}.
\]
\par
Thirdly, we try to find a vector field $\eta\in m_0\theta_0(2k-1)$
such that $W_2\frac{\partial}{\partial W_2}\ne \eta$ and $(W_2\circ
\varphi_k)\frac{\partial}{\partial W_2}\equiv \xi\circ \varphi_k$
(mod $T\mathcal {R}_e(\varphi_k))$.
\begin{eqnarray*}
(W_2\circ \varphi_k)\frac{\partial}{\partial W_2} & = &
\left(\sum_{i=1}^{k-1}v_iy^i\right)\frac{\partial}{\partial W_2} \\
{ } & \equiv & -\sum_{i=1}^{k-1}(V_i\circ
\varphi_k)\frac{\partial}{\partial V_i}\quad (\mbox{mod}\;
T\mathcal{R}_e(\varphi_k)).
\end{eqnarray*}
Thus, it follows that the following is a liftable vector field over
$\varphi_k$ where $\widetilde{\eta}_3=0$.
\[
\eta_3+\widetilde{\eta}_3=\sum_{i=1}^{k-1}V_i\frac{\partial}{\partial
V_i} +W_2\frac{\partial}{\partial W_2}.
\]
\par
Fourthly, since $y^i(W_1\circ \varphi_k)\frac{\partial}{\partial
W_1}\equiv -(W_1\circ \varphi_k)\frac{\partial}{\partial U_i}$ (mod
$T\mathcal{A}_e(\varphi_k)$) for any $i$ $(1\le i\le k-2)$, we try
to find a vector field $\eta_i\in m_0\theta_0(2k-1)$ such that
$-W_1\frac{\partial}{\partial U_i}\ne \eta_i$ and $y^i(W_1\circ
\varphi_k)\frac{\partial}{\partial W_1}\equiv \eta_i\circ \varphi_k$
(mod $T\mathcal
{R}_e(\varphi_k)+\varphi_k^2m_0^2\theta(\varphi_k)$).
\begin{eqnarray*}
{ } & { } & y^i(W_1\circ \varphi_k)\frac{\partial}{\partial W_1} \\
{ } & = &
y^i\left(y^k+\sum_{j=1}^{k-2}u_jy^j\right)\frac{\partial}{\partial W_1} \\
{ } & \equiv & \left(-\frac{1}{k}\sum_{j=1}^{k-2}u_jy^{i+j}+
\sum_{j=1}^{k-2}u_jy^{i+j}\right)\frac{\partial}{\partial W_1}
-\sum_{j=1}^{k-1}j v_j y^{i+j}\frac{\partial}{\partial W_2}
\;\;  (\mbox{mod}\; T\mathcal{R}_e(\varphi_k)).   \\
{ } & \equiv &
-\sum_{j=1}^{k-2-i}\frac{(k-j)}{k}(U_j\circ\varphi_k)\frac{\partial}{\partial
U_{i+j}}
+\sum_{j=1}^{k-1-i}j(V_j\circ\varphi_k)\frac{\partial}{\partial V_{i+j}} \\
{ } & { } & \qquad\qquad\qquad\qquad\qquad\qquad\qquad\qquad
 (\mbox{mod}\; T\mathcal{R}_e(\varphi_k)
+\varphi_k^2m_0^2\theta(\varphi_k)).
\end{eqnarray*}
Thus, it follows that for any $i$ $(1\le i\le k-2)$ there exists a
liftable vector field over $\varphi_k$ having the following form
where $\widetilde{\eta}_{3+i}\in m_0^2\theta_0(p)$.
\begin{eqnarray*}
{ } & { } & \eta_{3+i}+\widetilde{\eta}_{3+i} \\
{ } & = & W_1\frac{\partial}{\partial
U_i}-\sum_{j=1}^{k-2-i}\frac{(k-j)}{k}U_j\frac{\partial}{\partial
U_{i+j}} +\sum_{j=1}^{k-1-i}jV_j\frac{\partial}{\partial V_{i+j}} +
\mbox{ higher terms}.
\end{eqnarray*}
\par
Fifthly, since $y^i(W_2\circ \varphi_k)\frac{\partial}{\partial
W_1}\equiv -(W_2\circ \varphi_k)\frac{\partial}{\partial U_i}$ (mod
$T\mathcal{A}_e(\varphi_k)$) for any $i$ $(1\le i\le k-2)$, we try
to find a vector field $\eta_i\in m_0\theta_0(2k-1)$ such that
$-W_2\frac{\partial}{\partial U_i}\ne \eta_i$ and $y^i(W_2\circ
\varphi_k)\frac{\partial}{\partial W_1}\equiv \eta_i\circ \varphi_k$
(mod $T\mathcal
{R}_e(\varphi_k)+\varphi_k^2m_0^2\theta(\varphi_k)$).
\begin{eqnarray*}
{ } & { } & y^i(W_2\circ \varphi_k)\frac{\partial}{\partial W_1}  \\
{ } & = &
\left(\sum_{j=1}^{k-1}v_jy^{i+j}\right)\frac{\partial}{\partial W_1} \\
{ } & \equiv & \left\{
\begin{array}{ll}
-\sum_{j=1}^{k-2-i}(V_j\circ \varphi_k)\frac{\partial}{\partial U_{i+j}}& { } \\
\qquad\;\; (\mbox{mod}\; T\mathcal{R}_e(\varphi_k)
+\varphi_k^2m_0^2\theta(\varphi_k)) & (1\le i\le k-3),  \\
0   \qquad  (\mbox{mod}\; T\mathcal{R}_e(\varphi_k)
+\varphi_k^2m_0^2\theta(\varphi_k)) & (i=k-2).
\end{array}
\right.
\end{eqnarray*}
Thus, it follows that for any $i$ $(1\le i\le k-3)$ there exists a
liftable vector field over $\varphi_k$ having the following form
where $\widetilde{\eta}_{k+1+i}\in m_0^2\theta_0(p)$,
\[
\eta_{k+1+i}+\widetilde{\eta}_{k+1+i}= W_2\frac{\partial}{\partial
U_i} -\sum_{j=1}^{k-2-i}V_j\frac{\partial}{\partial U_{i+j}} +
\mbox{ higher terms};
\]
and there exists a liftable vector field over $\varphi_k$ having the
following form, where $\widetilde{\eta}_{2k-1}\in m_0^2\theta_0(p)$.
\[
\eta_{2k-1}+\widetilde{\eta}_{2k-1}=W_2\frac{\partial}{\partial
U_{k-2}}+ \mbox{ higher terms}.
\]
\par
Sixthly, since $y^i(W_2\circ \varphi_k)\frac{\partial}{\partial
W_2}\equiv -(W_2\circ \varphi_k)\frac{\partial}{\partial V_i}$ (mod
$T\mathcal{A}_e(\varphi_k)$) for any $i$ $(1\le i\le k-1)$, we try
to find a vector field $\eta_i\in m_0\theta_0(2k-1)$ such that
$-W_2\frac{\partial}{\partial V_i}\ne \eta_i$ and $y^i(W_2\circ
\varphi_k)\frac{\partial}{\partial W_2}\equiv \eta_i\circ \varphi_k$
(mod $T\mathcal
{R}_e(\varphi_k)+\varphi_k^2m_0^2\theta(\varphi_k)$).
\begin{eqnarray*}
{ } & { } & y^i(W_2\circ \varphi_k)\frac{\partial}{\partial W_2} \\
{ } & = &
\left(\sum_{j=1}^{k-1}v_jy^{i+j}\right)\frac{\partial}{\partial W_2} \\
{ } & \equiv & \left\{
\begin{array}{ll}
-\sum_{j=1}^{k-1-i}(V_j\circ \varphi_k)\frac{\partial}{\partial V_{i+j}}& { } \\
\qquad\;\; (\mbox{mod}\; T\mathcal{R}_e(\varphi_k)
+\varphi_k^2m_0^2\theta(\varphi_k)) & (1\le i\le k-2),  \\
0   \qquad  (\mbox{mod}\; T\mathcal{R}_e(\varphi_k)
+\varphi_k^2m_0^2\theta(\varphi_k)) & (i=k-1).
\end{array}
\right.
\end{eqnarray*}
Thus, it follows that for any $i$ $(1\le i\le k-2)$ there exists a
liftable vector field over $\varphi_k$ having the following form
where $\widetilde{\eta}_{2k-1+i}\in m_0^2\theta_0(p)$,
\[
\eta_{2k-1+i}+\widetilde{\eta}_{2k-1+i}= W_2\frac{\partial}{\partial
V_i} -\sum_{j=1}^{k-1-i}V_j\frac{\partial}{\partial V_{i+j}} +
\mbox{ higher terms};
\]
and there exists a liftable vector field over $\varphi_k$ having the
following form where $\widetilde{\eta}_{3k-2}\in m_0^2\theta_0(p)$.
\[
\eta_{3k-2}+\widetilde{\eta}_{3k-2}= W_2\frac{\partial}{\partial
V_{k-1}}+ \mbox{ higher terms}.
\]
\par
Finally, set
\begin{eqnarray*}
\Pi & = & \mathbb{K}W_1\frac{\partial}{\partial W_1}+
\mathbb{K}W_2\frac{\partial}{\partial W_1}+
\mathbb{K}W_2\frac{\partial}{\partial W_2} \\
{ } & { } & \qquad +
\sum_{i=1}^{k-2}\mathbb{K}W_1\frac{\partial}{\partial U_i}+
\sum_{i=1}^{k-2}\mathbb{K}W_2\frac{\partial}{\partial U_i}+
\sum_{i=1}^{k-1}\mathbb{K}W_2\frac{\partial}{\partial V_i}.
\end{eqnarray*}
Then, $\Pi$ is a $(3k-2)$-dimensional $\mathbb{K}$-vector space.
Let $\pi: \theta_0(p)\to \Pi$ be the canonical projection.   Then,
we see easily that $\pi(\eta_i+\widetilde{\eta}_i)$ $(1\le i\le
3k-2)$ constitute a basis of $\Pi$. Thus,
$\eta_i+\widetilde{\eta}_i$ $(1\le i\le 3k-2)$ constitute a set of
generators for the module of vector fields liftable over
$\varphi_k$.
\subsection{$Lift(\psi_n)$ for $\psi_n(v_1, \ldots, v_{n-1},y)=(v_1, \ldots, v_{n-1}, y^2,  v_1y, \ldots, v_{n-1}y)$}
\label{subsection 6.2}\qquad \\
We let $(V_1, \ldots, V_{n-1}, W, X_1, \ldots, X_{n-1})$ be the
standard coordinates of $\mathbb{K}^{2n-1}$. Since
${}_0\overline{\omega}\psi_n$ is bijective we first look for a basis
of $\mbox{ker}({}_1\overline{\omega}\psi_n)$. We can find out easily
a basis of $\mbox{ker}({}_1\overline{\omega}\psi_n)$ which is (for
instance) the following:
\begin{eqnarray*}
& { } & V_i\frac{\partial}{\partial V_j}+X_i\frac{\partial}{\partial
X_j}+m_0^2\theta_0(2n-1)
\qquad (1\le i, j\le n-1), \\
& { } &
X_i\frac{\partial}{\partial V_j}+m_0^2\theta_0(2n-1) \qquad (1\le i, j\le n-1), \\
& { } &
2X_i\frac{\partial}{\partial W}+m_0^2\theta_0(2n-1) \qquad (1\le i\le n-1), \\
& { } & 2W\frac{\partial}{\partial
W}+\sum_{j=1}^{n-1}X_j\frac{\partial}{\partial
X_j}+m_0^2\theta_0(2n-1).
\end{eqnarray*}
Since any component function of $\psi_n$ is a monomial, we can
determine easily the desired higher terms of liftable vector fields
and thus we see that the following constitute a set of generators
for the module of vector fields liftable over $\psi_n$.
\begin{eqnarray*}
& { } & V_i\frac{\partial}{\partial V_j}+X_i\frac{\partial}{\partial
X_j}
\qquad (1\le i, j\le n-1), \\
& { } &
X_i\frac{\partial}{\partial V_j}+V_iW\frac{\partial}{\partial X_j} \qquad (1\le i, j\le n-1), \\
& { } &
2X_i\frac{\partial}{\partial W}+\sum_{j=1}^{n-1}V_iV_j\frac{\partial}{\partial X_j} \qquad (1\le i\le n-1), \\
& { } & 2W\frac{\partial}{\partial
W}+\sum_{j=1}^{n-1}X_j\frac{\partial}{\partial X_j}.
\end{eqnarray*}
\subsection{$Lift(\phi)$ for
$\phi(x,y)=(x, y^2, y^3+x y)$
}
\label{subsection 6.3}\qquad \\
Let $(x,y)$, $(V,W,X)$ be the standard coordinates of $\mathbb{K}^2$
and $\mathbb{K}^3$ respectively, and let $\phi: (\mathbb{K}^2,0)\to
(\mathbb{K}^3,0)$ be the mono-germ defined by
\[
\phi(x,y)=(x, y^2, y^3+xy).
\]
Set,
\[
h(x,y)=(x+y^2, y) \quad \mbox{and}\quad H(V,W,X)=(V-W, W,X).
\]
Then, both $h: \textcolor{black}{\mathbb{K}^2\to \mathbb{K}^2}$ and
$H: \textcolor{black}{\mathbb{K}^3\to \mathbb{K}^3}$ are
\textcolor{black}{analytic} diffeomorphisms and preserve the origin.
Moreover, we have the following:
\[
\phi(x,y)=H\circ \psi_{\textcolor{black}{2}}\circ h(x,y),
\]
where $\psi_{\textcolor{black}{2}}$ is the mono-germ defined in
Example \ref{example 3.3}. By this equality, $f$ is
$\mathcal{A}$-equivalent to $\psi_{\textcolor{black}{2}}$. As same
as $\psi_{\textcolor{black}{2}}$, $f$ is often used as the normal
form of {\it Whitney umbrella} \textcolor{black}{from $\mathbb{K}^2$
to $\mathbb{K}^3$}. By Subsection \ref{subsection 6.2}, we have the
following:
\[
Lift(\psi_{\textcolor{black}{2}})=\left\langle
V\frac{\partial}{\partial V}+X\frac{\partial}{\partial X},
X\frac{\partial}{\partial V}+V W \frac{\partial}{\partial X},
 2X\frac{\partial}{\partial W}+V^2\frac{\partial}{\partial X},
2W\frac{\partial}{\partial W}+X\frac{\partial}{\partial X}
\right\rangle_{C_0}.
\]
Thus, by using the following lemma, $Lift(\phi)$ can be
characterized as the $C_0$-module generated by the following $4$
vector fields:
\begin{eqnarray*}
 { } & (V+W)\frac{\partial}{\partial V}+X\frac{\partial}{\partial X},
X\frac{\partial}{\partial V}+(V+W)W\frac{\partial}{\partial X}, \\
{ } &  -2X\frac{\partial }{\partial V}+2X\frac{\partial}{\partial
W}+(V+W)^2\frac{\partial}{\partial X}, -2W\frac{\partial}{\partial
V}+2W\frac{\partial}{\partial W}+X\frac{\partial}{\partial X}.
\end{eqnarray*}
\begin{lemma}\label{liftable lemma}
Let $S$ be a finite subset $\{s_1, \ldots, s_{r}\}$ $(s_i\ne
s_j\mbox{ if }i\ne j)$ and let $f=\textcolor{black}{\{}f_1, \ldots,
f_r\textcolor{black}{\}}, g=\textcolor{black}{\{}g_1, \ldots,
g_r\textcolor{black}{\}}: (\mathbb{K}^n, S)\to (\mathbb{K}^p,0)$ be
two \textcolor{black}{analytic} multigerms. Suppose that there exist
germs of analytic diffeomorphisms $h_i: (\mathbb{K}^n, s_i)\to
(\mathbb{K}^n, s_i)$ $(1\le i\le r)$ and $H: (\mathbb{K}^p, 0)\to
(\mathbb{K}^p, 0)$ such that $g=H\circ f\circ h$, where $h:
(\mathbb{K}^n, S)\to (\mathbb{K}^n, S)$ is the map-germ whose
restriction to $(\mathbb{K}^n, s_i)$ is $h_i$. Then, the mapping
$L_{(f,g)}: Lift(f)\to Lift(g) $ defined by $L_{(f,g)}(\eta)=dH\circ
\eta\circ H^{-1}$ is well-defined and bijective.
\end{lemma}
\par
\smallskip
\noindent \underline{\it Proof of Lemma \ref{liftable lemma}}\qquad
\par
Let $\eta$ be a liftable vector field over $f$.     By definition,
there exists $\xi\in\theta_S(n)$ such that $\eta\circ f = tf\circ
\xi$.   Since the equality $g=H\circ f\circ h$ holds, we have the
following:
\[
\eta\circ (H^{-1}\circ g\circ h^{-1})=t(H^{-1}\circ g\circ
h^{-1})\circ \xi.
\]
Hence, we have the following:
\[
(dH\circ \eta\circ H^{-1})\circ g=tg\circ (dh^{-1}\circ \xi\circ h).
\]
This shows that $(dH\circ \eta\circ H^{-1})$ is a liftable vector
field over $g$. Hence, the mapping $L_{(f,g)}$ is well-defined.
\par
Since injectivity of $L_{(f,g)}$ is clear, it is sufficient to show
that $L_{(f,g)}$ is surjective. Let $\widetilde{\eta}$
\textcolor{black}{be} a liftable vector field of $g$. The above
argument shows that $d(H^{-1})\circ \widetilde{\eta}\circ H$ is a
liftable vector field of $f$.   Since $L_{(f,g)}(d(H^{-1})\circ
\widetilde{\eta}\circ H)=\widetilde{\eta}$, it follows that
$L_{(f,g)}$ is surjective.
  \hfill $\Box$
\subsection{$Lift(f)$ for $f(x,y)=\{(x,y^2),(x^2,y)\}$}
\label{subsection 6.4}\qquad \\
Let $(X, Y)$ be the standard coordinates of $\mathbb{K}^{2}$. Since
${}_0\overline{\omega}f$ is bijective we first look for a basis of
$\mbox{ker}({}_1\overline{\omega}f)$. We can find out easily a basis
of $\mbox{ker}({}_1\overline{\omega}f)$ which is (for instance) the
following:
\[
X\frac{\partial}{\partial X}+m_0^2\theta_0(2), \;
Y\frac{\partial}{\partial Y}+m_0^2\theta_0(2).
\]
Since any component function of $f$ is a monomial, we can determine
easily the desired higher terms of liftable vector fields and thus
we see that the following constitute a set of generators for the
module of vector fields liftable over $f$.
\[
X\frac{\partial}{\partial X}, \; Y\frac{\partial}{\partial Y}.
\]
\subsection{$Lift(f)$ for $f(x)=\{(x^2,x^3),(x^3,x^2)\}$}
\label{subsection 6.5}\qquad \\
Recall that the multigerm $f$ of Example \ref{example 3.5}.2 is
$f_1(x)=(x^2,x^3), f_2(x)= (x^3,x^2)$. Let $(X, Y)$ be the standard
coordinates of $\mathbb{K}^{2}$. Since ${}_1\overline{\omega}f$ is
bijective we first look for a basis of
$\mbox{ker}({}_2\overline{\omega}f)$. We can find out easily a basis
of $\mbox{ker}({}_2\overline{\omega}f)$ which is (for instance) the
following:
\[
 6XY\frac{\partial}{\partial X}+4Y^2\frac{\partial}{\partial Y}+m_0^3\theta_0(2),
 \;
4X^2\frac{\partial}{\partial X}+6XY\frac{\partial}{\partial
Y}+m_0^3\theta_0(2).
\]
\par
Set $\xi_{1,1,1}=3x^4\frac{\partial}{\partial x}$,
$\xi_{1,2,1}=2x^3\frac{\partial}{\partial x}$ and
$\eta_{1,1}=6XY\frac{\partial}{\partial
X}+4Y^2\frac{\partial}{\partial Y}$. Then, we have the following:
\begin{eqnarray*}
{ } & { } & \eta_{1,1}\circ f_1 - df_1\circ \xi_{1,1,1}=-5x^6\frac{\partial}{\partial Y},  \\
{ } & { } & \eta_{1,1}\circ f_2 - df_2\circ \xi_{1,2,1}=0.
\end{eqnarray*}
Set $\eta_{1,2}=5X^3\frac{\partial}{\partial Y}$.   Then we have the
following:
$$
 (\eta_{1,1}+\eta_{1,2})\circ f_1 - df_1\circ \xi_{1,1,1}=0,  \eqno{(6.1)}
$$
$$
(\eta_{1,1}+\eta_{1,2})\circ f_2 - df_2\circ
\xi_{1,2,1}=5x^9\frac{\partial}{\partial Y}. \eqno{(6.2)}
$$
Set $\eta_{1,3}=-5XY^3\frac{\partial}{\partial Y}$ and
$\xi_{1,1,2}=-\frac{5}{3}x^9\frac{\partial}{\partial x}$. Then we
have the following:
 \begin{eqnarray*}
{ } & { } & (\eta_{1,1}+\eta_{1,2}+\eta_{1,3})\circ f_1 - df_1\circ
(\xi_{1,1,1}+\xi_{1,1,2})
=\frac{10}{3}x^{10}\frac{\partial}{\partial X},  \\
{ } & { } & (\eta_{1,1}+\eta_{1,2}+\eta_{1,3})\circ f_2 - df_2\circ
\xi_{1,2,1}=0.
\end{eqnarray*}
Set $\eta_{1,4}=-\frac{10}{3}X^2Y^2\frac{\partial}{\partial X}$ and
$\xi_{1,2,2}=-\frac{10}{9}x^8\frac{\partial}{\partial x}$. Then we
have the following:
$$
(\eta_{1,1}+\eta_{1,2}+\eta_{1,3}+\eta_{1,4})\circ f_1 - df_1\circ
(\xi_{1,1,1}+\xi_{1,1,2})=0,
\eqno{(6.3)} \\
$$
$$
 (\eta_{1,1}+\eta_{1,2}+\eta_{1,3}+\eta_{1,4})\circ f_2 - df_2\circ (\xi_{1,2,1}+\xi_{1,2,2})=\frac{20}{9}x^9\frac{\partial}{\partial Y}.
\eqno{(6.4)}
$$
Note that the right hand side of (5.3) (resp., the right hand side
of (5.4)) is the right hand side of (5.1) (resp., the right hand
side of (5.2)) multiplied by $\left(\frac{2}{3}\right)^2$. Thus, the
following vector field $\eta_1$ must be liftable over $f$.
\begin{eqnarray*}
\eta_{1} & = &
\eta_{1,1}+\eta_{1,2}+\left(1+\left(\frac{2}{3}\right)^2+\left(\frac{2}{3}\right)^4+\cdots\right)\left(\eta_{1,3}+\eta_{1,4}\right) \\
{ } & = & (6XY-6X^2Y^2)\frac{\partial}{\partial X}+
(4Y^2+5X^3-9XY^3)\frac{\partial}{\partial Y}.
\end{eqnarray*}
\par
Next, Set $\xi_{2,1,1}=2x^3\frac{\partial}{\partial x}$,
$\xi_{2,2,1}=3x^4\frac{\partial}{\partial x}$ and
$\eta_{2,1}=4X^2\frac{\partial}{\partial
X}+6XY\frac{\partial}{\partial Y}$. Then, we have the following:
\begin{eqnarray*}
{ } & { } & \eta_{2,1}\circ f_1 - df_1\circ \xi_{2,1,1}=0,  \\
{ } & { } & \eta_{2,1}\circ f_2 - df_2\circ
\xi_{2,2,1}=-5x^6\frac{\partial}{\partial X}.
\end{eqnarray*}
Set $\eta_{2,2}=5Y^3\frac{\partial}{\partial X}$.   Then we have the
following:
$$
 (\eta_{2,1}+\eta_{2,2})\circ f_1 - df_1\circ \xi_{2,1,1}=5x^9\frac{\partial}{\partial X},
\eqno{(6.5)}\\
$$
$$
 (\eta_{2,1}+\eta_{2,2})\circ f_2 - df_2\circ \xi_{2,2,1}=0.
\eqno{(6.6)}
$$
Set $\eta_{2,3}=-5X^3Y\frac{\partial}{\partial X}$ and
$\xi_{2,2,2}=-\frac{5}{3}x^9\frac{\partial}{\partial x}$. Then we
have the following:
 \begin{eqnarray*}
{ } & { } & (\eta_{2,1}+\eta_{2,2}+\eta_{2,3})\circ f_1 - df_1\circ (\xi_{2,1,1})=0,  \\
{ } & { } & (\eta_{2,1}+\eta_{2,2}+\eta_{2,3})\circ f_2 - df_2\circ
(\xi_{2,2,1}+\xi_{2,2,2})=\frac{10}{3}x^{10}\frac{\partial}{\partial
Y}.
\end{eqnarray*}
Set $\eta_{2,4}=-\frac{10}{3}X^2Y^2\frac{\partial}{\partial Y}$ and
$\xi_{2,1,2}=-\frac{10}{9}x^8\frac{\partial}{\partial x}$. Then we
have the following:
$$
 (\eta_{2,1}+\eta_{2,2}+\eta_{2,3}+\eta_{2,4})\circ f_1 - df_1\circ (\xi_{2,1,1}+\xi_{2,1,2})=\frac{20}{9}x^9\frac{\partial}{\partial X},
\eqno{(6.7)}\\
$$
$$
 (\eta_{2,1}+\eta_{2,2}+\eta_{2,3}+\eta_{2,4})\circ f_2 - df_2\circ (\xi_{2,2,1}+\xi_{2,2,2})=0.
\eqno{(6.8)}
$$
Note that the right hand side of (5.7) (resp., the right hand side
of (5.8)) is the right hand side of (5.5) (resp., the right hand
side of (5.6)) multiplied by $\left(\frac{2}{3}\right)^2$. Thus, the
following vector field $\eta_2$ must be liftable over $f$.
\begin{eqnarray*}
\eta_{2} & = &
\eta_{2,1}+\eta_{2,2}+\left(1+\left(\frac{2}{3}\right)^2+\left(\frac{2}{3}\right)^4+\cdots\right)\left(\eta_{2,3}+\eta_{2,4}\right) \\
{ } & = & (4X^2+5Y^3-9X^3Y)\frac{\partial}{\partial X}+
(6XY-6X^2Y^2)\frac{\partial}{\partial Y}.
\end{eqnarray*}
\par
Therefore, the following constitute a set of generators for the
module of vector fields liftable over $f$.
\begin{eqnarray*}
\eta_1 & = & (6XY-6X^2Y^2)\frac{\partial}{\partial X}+
(4Y^2+5X^3-9XY^3)\frac{\partial}{\partial Y}, \\
\eta_2 & = & (4X^2+5Y^3-9X^3Y)\frac{\partial}{\partial X}+
(6XY-6X^2Y^2)\frac{\partial}{\partial Y}.
\end{eqnarray*}
\subsection{$Lift(f)$ for $f(y)= (y^2,0)$}
\label{subsection 6.6}\qquad \\
Let $f: (\mathbb{K},0)\to (\mathbb{K}^2, 0)$ be the mono-germ
defined by $f(y)=(y^2,0)$. As an application of Theorem \ref{theorem
2}, we obtain all liftable vector fields over $f$.
\par
Let $(Y,U)$ be the standard coordinates of the target space of $f$.
It is easy to see the following:
\[
\theta_S(f)=T\mathcal{K}_e(f)+\mathbb{K}^2+y\frac{\partial}{\partial
U}
\]
Since
$\dim_{\mathbb{K}}\theta_S(f)/\left(T\mathcal{K}_e(f)+\mathbb{K}^2\right)=1$,
by Mather's constructing method of stable mono-germs
(\cite{mather4}), the mono-germ $F(x, y)=(x, y^2, xy)$ is a
one-parameter stable unfolding of $f$. Notice that $F$ is exactly
the same as the mono-germ $\psi_2$ defined in Subsection
\ref{subsection 6.2}. Let $(X,Y,U)$ be the standard coordinates of
the target space of $F$.  Let $g: (\mathbb{K}\times\mathbb{K}^2,
(0,0))\to (\mathbb{K}\times\mathbb{K}^2, (0,0))$ be defined by $g(x,
y,u)=(x^2, y,u)$. Then, $Lift(g)=\langle X\frac{\partial}{\partial
X}, \frac{\partial}{\partial Y}, \frac{\partial}{\partial
U}\rangle_{C_0}$. Set $\widetilde{\eta}_1=X\frac{\partial}{\partial
X}+U\frac{\partial}{\partial U}$,
$\widetilde{\eta}_2=U\frac{\partial}{\partial X}+X Y
\frac{\partial}{\partial U}$, $\widetilde{\eta}_3=
2U\frac{\partial}{\partial Y}+X^2\frac{\partial}{\partial U}$ and
$\widetilde{\eta}_4=2Y\frac{\partial}{\partial
Y}+U\frac{\partial}{\partial U}$. By Subsection \ref{subsection
6.2}, we have the following:
\begin{eqnarray*}
Lift(F) & = & \left\langle \widetilde{\eta}_1, \widetilde{\eta}_2,
\widetilde{\eta}_3, \widetilde{\eta}_4
\right\rangle_{C_0} \\
{ } & = & \left\{
\widetilde{\alpha}_1\left(X\frac{\partial}{\partial
X}+U\frac{\partial}{\partial U}\right)
+ \widetilde{\alpha}_2\left(U\frac{\partial}{\partial X}+X Y\frac{\partial}{\partial U}\right) \right.\\
{ } & { } & \quad \left. +
\widetilde{\alpha}_3\left(2U\frac{\partial}{\partial
Y}+X^2\frac{\partial}{\partial U}\right) +
\widetilde{\alpha_4}\left(2Y\frac{\partial}{\partial
Y}+U\frac{\partial}{\partial U}\right) \right\},
\end{eqnarray*}
where $\widetilde{\alpha}_i$ $(1\le i\le 4)$ are
\textcolor{black}{analytic}  function-germs of three variables $X,
Y, U$. Thus, we have the following:
\begin{eqnarray*}
{ } & { } & Lift(F)\cap Lift(g) \\
{ } & = & \left\{
\widetilde{\alpha}_1\left(X\frac{\partial}{\partial
X}+U\frac{\partial}{\partial U}\right) +
\widetilde{\alpha}_2\left(U\frac{\partial}{\partial X} +X
Y\frac{\partial}{\partial U}\right) +
\widetilde{\alpha}_3\left(2U\frac{\partial}{\partial
Y}+X^2\frac{\partial}{\partial U}\right)
\right.\\
{ } & { } & \quad \left.\left. +
\widetilde{\alpha}_4\left(2Y\frac{\partial}{\partial Y}
+U\frac{\partial}{\partial U}\right) \right|
\widetilde{\alpha}_1X+\widetilde{\alpha}_2U \mbox{ can be divided by
}X
\right\} \\
{ } & = & \left\{
\widetilde{\alpha}_1\left(X\frac{\partial}{\partial
X}+U\frac{\partial}{\partial U}\right) +
\widetilde{\alpha}_2\left(U\frac{\partial}{\partial X} +X
Y\frac{\partial}{\partial U}\right) +
\widetilde{\alpha}_3\left(2U\frac{\partial}{\partial Y}
+X^2\frac{\partial}{\partial U}\right)
\right. \\
{ } & { } & \quad \left. \left. +
\widetilde{\alpha}_4\left(2Y\frac{\partial}{\partial Y}
+U\frac{\partial}{\partial U}\right) \right|
\widetilde{\alpha}_2 \mbox{ can be divided by }X \right\}.
\end{eqnarray*}
Define $\alpha_i: (\mathbb{K}^2, 0)\to \mathbb{K}$ $(1\le i\le 4)$
by $\alpha_i(Y,U)=\widetilde{\alpha}_i(0, Y, U)$. Then, by Theorem
\ref{theorem 2}, we have the following:
\begin{eqnarray*}
{ } & { } & Lift(f) \\
{ } & = & \left\{ \alpha_1U\frac{\partial}{\partial U} +
2\alpha_3U\frac{\partial}{\partial Y} +
\alpha_4\left(2Y\frac{\partial}{\partial
Y}+U\frac{\partial}{\partial U}\right)
\right\} \\
{ } & = & \left\langle U\frac{\partial}{\partial U},
U\frac{\partial}{\partial Y}, Y\frac{\partial}{\partial Y}
\right\rangle_{C_0}.
\end{eqnarray*}
Since three vector fields $U\frac{\partial}{\partial U},
U\frac{\partial}{\partial Y}, Y\frac{\partial}{\partial Y}$ are
linearly independent, the minimal number of generators for $Lift(f)$
is $3$, which is strictly greater than the dimension of the target
space of $f$.
\subsection{$Lift(f_k)$ for $f_k(y)= (y^2,y^{2k+1})$ $(k\ge 1)$
}
\label{subsection 6.7}\qquad \\
Let $f_k: (\mathbb{K},0)\to (\mathbb{K}^2, 0)$ $(k\ge 1)$ be the
mono-germ defined by $f_k(y)=(y^2,y^{2k+1})$. As an application of
Theorem \ref{theorem 2}, we obtain all liftable vector fields over
$f_k$.
\par
Let $(Y,U)$ be the standard coordinates of the target space of
$f_k$. It is easy to see the following:
\[
\theta_S(f_k)=T\mathcal{K}_e(f_k)+\mathbb{K}^2+y\frac{\partial}{\partial
U}
\]
Since
$\dim_{\mathbb{K}}\theta_S(f_k)/\left(T\mathcal{K}_e(f_k)+\mathbb{K}^2\right)=1$,
by Mather's constructing method of stable mono-germs
(\cite{mather4}), the mono-germ $F_k(x, y)=(x, y^2, y^{2k+1}+xy)$ is
a one-parameter stable unfolding of $f_k$. Set $F(x,y)=(x,y^2,xy)$.
Let $(X,Y,U)$ be the standard coordinates of the target space of
$F_k$. Set $h_k(x,y)=(x+y^{2k}, y)$ and $H_k(X,Y,U)=(X-Y^k, Y, U)$.
Then, both $h_k$ and $H_k$ are \textcolor{black}{analytic}
diffeomorphisms preserving the origin, and we have that
$F_k=H_k\circ F\circ h_k$. Set
$\widetilde{\eta}_1=X\frac{\partial}{\partial
X}+U\frac{\partial}{\partial U}$,
$\widetilde{\eta}_2=U\frac{\partial}{\partial X}+ X Y
\frac{\partial}{\partial U}$, $\widetilde{\eta}_3=
2U\frac{\partial}{\partial Y}+X^2\frac{\partial}{\partial U}$ and
$\widetilde{\eta}_4=2Y\frac{\partial}{\partial
Y}+U\frac{\partial}{\partial U}$. By Subsection \ref{subsection
6.2}, we have the following:
\[
Lift(F)=\langle \widetilde{\eta}_1, \widetilde{\eta}_2,
\widetilde{\eta}_3, \widetilde{\eta}_4 \rangle_{C_0}.
\]
Set $\eta_1=d H_k\circ \widetilde{\eta}_1\circ H_k^{-1}$, $\eta_2=d
H_k\circ \widetilde{\eta}_2\circ H_k^{-1}$, $\eta_3=d H_k\circ
\widetilde{\eta}_3\circ H_k^{-1}$ and $\eta_4=d H_k\circ
\widetilde{\eta}_4\circ H_k^{-1}$. By Lemma \ref{liftable lemma}, we
have the following:
\[
Lift(F_k)=\langle \eta_1, \eta_2, \eta_3, \eta_4 \rangle_{C_0}.
\]
By calculations, we have the following:
\[
\left\{
\begin{array}{ccl}
\eta_1(X,Y,U) & = &
(X+Y^k)\frac{\partial}{\partial X} + U\frac{\partial}{\partial U}, \\
\eta_2(X,Y,U) & = &
U\frac{\partial}{\partial X} + (X+Y^k)Y\frac{\partial}{\partial U},  \\
\eta_3(X,Y,U) & = & -2kY^{k-1}U\frac{\partial}{\partial X} +
2U\frac{\partial}{\partial Y}+(X+Y^k)^2\frac{\partial}{\partial U}, \\
\eta_4(X,Y,U) & = & -2kY^k\frac{\partial}{\partial X} +
2Y\frac{\partial}{\partial Y}+U\frac{\partial}{\partial U}.
\end{array}
\right.
\]
\par
Let $g: (\mathbb{K}\times\mathbb{K}^2, (0,0))\to
(\mathbb{K}\times\mathbb{K}^2, (0,0))$ be defined by $g(x,
y,u)=(x^2, y,u)$. Then, $Lift(g)=\langle X\frac{\partial}{\partial
X}, \frac{\partial}{\partial Y}, \frac{\partial}{\partial
U}\rangle_{C_0}$. Thus, we have the following:
\begin{eqnarray*}
{ } & { } & Lift(F_k)\cap Lift(g)  \\
{ } & = & \left\{ \left.\sum_{i=1}^4 \widetilde{\alpha}_i\eta_i\;
\right|\; \widetilde{\alpha}_1(X+Y^k)+\widetilde{\alpha}_2U
-2k\widetilde{\alpha}_3Y^{k-1}U-2k\widetilde{\alpha}_4Y^k \mbox{ can
be divided by }X \right\},
\end{eqnarray*}
where $\widetilde{\alpha}_i$ $(1\le i\le 4)$ are
\textcolor{black}{analytic}  function-germs of three variables $X,
Y, U$. Define $\alpha_i: (\mathbb{K}^2, 0)\to \mathbb{K}$ $(1\le
i\le 4)$ by $\alpha_i(Y,U)=\widetilde{\alpha}_i(0, Y, U)$. Then, by
Theorem \ref{theorem 2}, $Lift(\textcolor{black}{f_k})$ can be
characterized as follows:
\begin{eqnarray*}
{ } & { } & Lift(\textcolor{black}{f_k}) \\
{ } & = &
\left\{ 2(U\alpha_3+Y\alpha_4)\frac{\partial}{\partial Y}\right. \\
{ } & { } & \qquad \left.\left. +
(U(\alpha_1+\alpha_4)+Y^{2k}\alpha_3+Y^{k+1}\alpha_2)\frac{\partial}{\partial
U} \right|
Y^k\alpha_1+U\alpha_2-2k Y^{k-1}U\alpha_3-2kY^k\alpha_4=0
\right\} \\
{ } & = &
\left\{ 2(U\alpha_3+Y\alpha_4)\frac{\partial}{\partial Y} \right. \\
{ } & { } & \qquad \left.\left. +
(U(\alpha_1+\alpha_4)+Y^{2k}\alpha_3+Y^{k+1}\alpha_2)\frac{\partial}{\partial
U} \right|
Y^k(\alpha_1-2k\alpha_4)+U(\alpha_2-2kY^{k-1}\alpha_3)=0
\right\} \\
{ } & = &
\left\{2(U\alpha_3+Y\alpha_4)\frac{\partial}{\partial Y}\right. \\
{ } & { } & \qquad \left.\left. +
(U(\alpha_1+\alpha_4)+Y^{2k}\alpha_3+Y^{k+1}\alpha_2)\frac{\partial}{\partial
U} \right|
\alpha_1-2k\alpha_4=U\beta,\; \alpha_2-2k Y^{k-1}\alpha_3=-Y^k\beta
\right\}\\
{ } & = & \left\{2(U\alpha_3+Y\alpha_4)\frac{\partial}{\partial Y}+
\left(U(1+2k)\alpha_4+U^2\beta)+Y^{2k}(1+2k)\alpha_3-Y^{2k+1}\beta\right)\frac{\partial}{\partial
U}
\right\}  \\
{ } & = &
\left\{\alpha_3\left(2U\frac{\partial}{\partial Y}+(1+2k)Y^{2k}\frac{\partial}{\partial U}\right) \right. \\
{ } & { } & \qquad\qquad \left. +
\alpha_4\left(2Y\frac{\partial}{\partial
Y}+(1+2k)U\frac{\partial}{\partial U}\right) +
\beta\left(U^2-Y^{2k+1}\right)\frac{\partial}{\partial U}
\right\}  \\
{ } & = & \left\{\left(\alpha_3-\frac{1}{1+2k}Y\beta\right)
\left(2U\frac{\partial}{\partial
Y}+(1+2k)Y^{2k}\frac{\partial}{\partial U}\right)
\right. \\
{ } & { } & \qquad\qquad \left. +
\left(\alpha_4+\frac{1}{1+2k}U\beta\right)
\left(2Y\frac{\partial}{\partial Y}+(1+2k)U\frac{\partial}{\partial
U}\right)
\right\}  \\
{ } & = & \left\langle 2U\frac{\partial}{\partial
Y}+(1+2k)Y^{2k}\frac{\partial}{\partial U}, \;
2Y\frac{\partial}{\partial Y}+(1+2k)U\frac{\partial}{\partial U}
\right\rangle_{C_0},
%
\end{eqnarray*}
where $\beta$ is an \textcolor{black}{analytic}  function-germ of
two variables $Y, U$.
\subsection{$Lift(S_k^{\pm})$ for $S_k^{\pm}(x,y)= (x, y^2, y^3\pm x^{k+1}y)$
$(k\ge 0)$
}
\label{subsection 6.8} \qquad \\
Let $S_k^{\pm}: (\mathbb{K}^2,0)\to (\mathbb{K}^3, 0)$ be the
mono-germ defined by $S_k^{\pm}(x,y)=(x,y^2,y^3\pm x^{k+1}y)$ $(k\ge
0)$. The mono-germ $S_k^{\pm}$ can be found in the classification
list of $\mathcal{A}$-simple mono-germ from the plane to $3$-space
due to Mond (\cite{mond}). Here, a multigerm $f: (\mathbb{K}^n,
S)\to (\mathbb{K}^p,0)$ is said to be $\mathcal{A}$-{\it simple}
if there exists a finite number of $\mathcal{A}$-equivalence classes
such that for any positive integer $d$ and any
\textcolor{black}{analytic} mapping $F: U\to V$ where $U\subset
\mathbb{K}^n\times \mathbb{K}^d$ is a neighbourhood of $S\times 0$,
$V\subset \mathbb{K}^p\times \mathbb{K}^d$ is a neighbourhood of
$(0,0)$, $F(x, \lambda)=(f_\lambda(x),\lambda)$ and the germ of
$f_0$ at $S$ is $f$, there exists a sufficiently small neighbourhood
$W_i\subset U$ of $(s_i,0)$ $(1\le i\le |S|)$ such that for
every $\{(x_1,\lambda), \cdots, (x_r,\lambda)\}$ $(r\le |S|)$ with
$(x_i, \lambda)\in W_i$ and $F(x_1, \lambda)=\cdots =F(x_r,
\lambda)$
the multigerm $f_\lambda : (\mathbb{K}^n,\{x_1,\cdots, x_r\})\to
(\mathbb{K}^p,f_\lambda(x_i))$ lies in one of these finite
$\mathcal{A}$-equivalence classes. As an application of Theorem
\ref{theorem 2}, we obtain all liftable vector fields over
$S_k^{\pm}$.
\par
Let $(X,Y,U)$ be the standard coordinates of the target space of
$S_k^{\pm}$. It is easy to see the following:
\[
\theta_S(S_k^{\pm})=T\mathcal{K}_e(S_k^{\pm})+\mathbb{K}^3+y\frac{\partial}{\partial
U}
\]
Since
$\dim_{\mathbb{K}}\left(\theta_S(S_k^{\pm})/T\mathcal{K}_e(S_k^{\pm})+\mathbb{K}^3)\right)\le
1$, by Mather's constructing method of stable mono-germs
(\cite{mather4}), the mono-germ $F_k^{\pm}(x, y, u)=(x, y^2, y^3\pm
x^{k+1}y+u y, u)$ is a one-parameter stable unfolding of
$S_k^{\pm}$.
Set $F(x,y,u)=(x,y^2,y^3+u y,u)$.
Let $(X,Y,U,V)$ be the standard coordinates of the target space of
$F_k$. Set $h_k^\pm(x,y,u)=(x, y, u\pm x^{k+1})$ and
$H_k^\pm(X,Y,U,V)=(X, Y, U, V\mp X^{k+1})$. Then, both $h_k^\pm$ and
$H_k^\pm$ are \textcolor{black}{analytic} diffeomorphisms preserving
the origin, and we have that $F_k^\pm=H_k^\pm\circ F\circ h_k^\pm$.
Set \textcolor{black} {
$\widetilde{\eta}_1=U\frac{\partial}{\partial
U}+(Y+V)\frac{\partial}{\partial V}$,
$\widetilde{\eta}_2=(Y+V)Y\frac{\partial}{\partial U}+ U
\frac{\partial}{\partial V}$, $\widetilde{\eta}_3=
2U\frac{\partial}{\partial Y} +(Y+V)^2\frac{\partial}{\partial U}
\textcolor{black}{-}2U\frac{\partial}{\partial V}$, $\widetilde{\eta}_4=
2Y\frac{\partial}{\partial Y}+U\frac{\partial}{\partial U}
-2Y\frac{\partial}{\partial V}$ and
 $\widetilde{\eta}_5=\frac{\partial}{\partial X}$.
} By calculations, we have the following:
\[
Lift(F)=\langle \widetilde{\eta}_1, \widetilde{\eta}_2,
\widetilde{\eta}_3, \widetilde{\eta}_4,
 \widetilde{\eta}_5
\rangle_{C_0}.
\]
Set $\eta_1=d \textcolor{black}{(}H_k^{\textcolor{black}{\pm}}\textcolor{black}{)}
\circ \widetilde{\eta}_1\circ 
\textcolor{black}{(}H_k^{\textcolor{black}{\pm}}\textcolor{black}{)}^{-1}$, 
$\eta_2=d \textcolor{black}{(}H_k^{\textcolor{black}{\pm}}\textcolor{black}{)}
\circ \widetilde{\eta}_2\circ 
\textcolor{black}{(}H_k^{\textcolor{black}{\pm}}\textcolor{black}{)}^{-1}$, 
$\eta_3=d \textcolor{black}{(}H_k^{\textcolor{black}{\pm}}\textcolor{black}{)}
\circ \widetilde{\eta}_3\circ 
\textcolor{black}{(}H_k^{\textcolor{black}{\pm}}\textcolor{black}{)}^{-1}$, 
$\eta_4=d \textcolor{black}{(}H_k^{\textcolor{black}{\pm}}\textcolor{black}{)}
\circ \widetilde{\eta}_4\circ 
\textcolor{black}{(}H_k^{\textcolor{black}{\pm}}\textcolor{black}{)}^{-1}$  
and 
$\eta_{\textcolor{black}{5}}=
d \textcolor{black}{(}H_k^{\textcolor{black}{\pm}}\textcolor{black}{)}
\circ \widetilde{\eta}_5\circ 
\textcolor{black}{(}H_k^{\textcolor{black}{\pm}}\textcolor{black}{)}^{-1}$.    
By Lemma \ref{liftable lemma}, we
have the following:
\[
Lift(F_k)=\langle \eta_1, \eta_2, \eta_3, \eta_4, \eta_5
\rangle_{C_0}.
\]
By calculations, we have the following: \textcolor{black} {
\[
\left\{
\begin{array}{ccl}
\eta_1(X,Y,U,V) & = &
U\frac{\partial}{\partial U} + \left(Y+V\pm X^{k+1}\right)\frac{\partial}{\partial V}, \\
\eta_2(X,Y,U,V) & = & Y(Y+V\pm X^{k+1})\frac{\partial}{\partial Y}
+U\frac{\partial}{\partial V},  \\
\eta_3(X,Y,U, V) & = & 2U\frac{\partial}{\partial Y} +\left(Y+V\pm
X^{k+1}\right)^2\frac{\partial}{\partial U}
\textcolor{black}{-} 2U\frac{\partial}{\partial V}, \\
\eta_4(X,Y,U, V) & = & 2Y\frac{\partial}{\partial Y} +
U\frac{\partial}{\partial U}-2Y\frac{\partial}{\partial V}, \\
\eta_5(X,Y,U, V) & = & \frac{\partial}{\partial X} \mp
(k+1)X^{k}\frac{\partial}{\partial V}.
\end{array}
\right.
\]
}
\par
Let $g: (\mathbb{K}^3\times\mathbb{K}, (0,0))\to (\mathbb{K}^3\times
\mathbb{K}, (0,0))$ be defined by $g(x, y,u,v)=(x, y,u, v^2)$. Then,
$Lift(g)=\langle \frac{\partial}{\partial X},
\frac{\partial}{\partial Y}, \frac{\partial}{\partial U},
V\frac{\partial}{\partial V}\rangle_{C_0}$. Thus, we have the
following: \textcolor{black} {
\[
Lift(F_k)\cap Lift(g)   =  \left\{
\left.\sum_{i=1}^5\widetilde{\alpha}_i\eta_i\; \right|\;
\Phi(X,Y,U,V) \mbox{ can be divided by } V \right\},
\]
} where $\widetilde{\alpha}_i$ $(1\le i\le 5)$ are
\textcolor{black}{analytic}  function-germs of four variables $X, Y,
U, V$ and $\Phi(X,Y,U,V)$ is given as follows: \textcolor{black} {
\begin{eqnarray*}
\Phi(X,Y,U,V) & = &
\widetilde{\alpha}_1(X,Y,U,V)\left(Y+V\pm X^{k+1}\right) \\
{ } & { } & \qquad +\widetilde{\alpha}_2(X,Y,U,V)U
\textcolor{black}{-}2\widetilde{\alpha}_3(X,Y,U,V)U \\
{ } & { } & \qquad
-\textcolor{black}{2}
\widetilde{\alpha}_4(X,Y,U,V)Y\mp(k+1)\widetilde{\alpha}_5(X,Y,U,V)X^k.
\end{eqnarray*}
} Define $\alpha_i: (\mathbb{K}^3, 0)\to \mathbb{K}$ $(1\le i\le 5)$
by $\alpha_i(X, Y,U)=\widetilde{\alpha}_i(X, Y, U, 0)$. Then, by
Theorem \ref{theorem 2}, any element of $Lift(S_k^{\pm})$ has the
following form: \textcolor{black} {
\begin{eqnarray*}
\alpha_1U\frac{\partial }{\partial U} + \alpha_2Y\left(Y\pm
X^{k+1}\right)\frac{\partial }{\partial Y}+
\alpha_3\left(2U\frac{\partial}{\partial Y}
+\left(Y\pm X^{k+1}\right)^2\frac{\partial }{\partial U}\right)  \\
\qquad +\alpha_4\left(2Y\frac{\partial }{\partial Y}
+U\frac{\partial }{\partial U} \right) +\alpha_5\frac{\partial
}{\partial X}.
\end{eqnarray*}
} And, by Theorem \ref{theorem 2} again, the unique restriction on
$\alpha_i$ $(1\le i\le 5)$ is as follows. \textcolor{black} {
\begin{condition}\label{condition 1}
\[
\Phi(X, Y, U, 0) = \alpha_1\left(Y\pm X^{k+1}\right) +\alpha_2U
\textcolor{black}{-}2\alpha_3U-2\alpha_4 Y \mp (k+1)\alpha_5 X^k = 0.
\]
\end{condition}
} 
In the case $k=0$, by Condition \ref{condition 1} it follows that
$\alpha_5$ can be expressed by using $\alpha_i$ $(1\le i\le 4)$.
Thus, it is easy to obtain four vector fields which constitute a
generators of $Lift(S_0^{\pm})$. \textcolor{black} { In the case
$k\ge 1$, by Condition \ref{condition 1}, we have the following
expressions:
\begin{eqnarray*}
\alpha_1-2\alpha_4 & = & \beta_1 X^k +\gamma U, \\
\alpha_2\textcolor{black}{-}2\alpha_3 & = & \beta_2 X^k -\gamma Y, \\
\end{eqnarray*}
where $\beta_1, \beta_2, \gamma$ are analytic function-germs
$(\mathbb{K}^3,0)\to \mathbb{K}$. Therefore, we have the following:
\[
\alpha_5 = \pm \frac{1}{k+1} \left( \beta_1Y+\beta_2U \pm
\left(2\alpha_4+\beta_1X^k+\gamma U\right)X \right).
\]
Hence, $Lift(S_k^{\pm})$ in the case $k\ge 1$ can be characterized
as follows:
\begin{eqnarray*}
{ } & { } & Lift(S_k^{\pm}) \\
{ } & = & \left\{ \left(2\alpha_4+\beta_1X^k+\gamma U\right)U
\frac{\partial}{\partial U} +\left( 
2\alpha_3+\beta_2X^k-\gamma
Y\right) \left(Y\pm X^{k+1}\right)Y \frac{\partial}{\partial U}
\right. \\
{ } & { } & \qquad + \alpha_3 \left(2U\frac{\partial}{\partial
Y}+\left(Y\pm X^{k+1}\right)^2\frac{\partial}{\partial U}\right) +
\alpha_4
\left(2Y\frac{\partial}{\partial Y}+U\frac{\partial}{\partial U}\right)  \\
{ } & { } & \qquad \left.\pm \frac{1}{k+1}\left(\beta_1 Y+\beta_2
U\pm \left( 2\alpha_4+\beta_1 X^k+\gamma U\right)X\right)
\frac{\partial}{\partial X}
\right\} \\
{ } & = & 
\left\langle 2U\frac{\partial}{\partial Y} +
\left(\textcolor{black}{3Y^2\pm 4X^{k+1}Y+X^{2k+2}} \right) \frac{\partial}{\partial U}, 
\; 
\frac{2X}{k+1}\frac{\partial}{\partial X}
+2Y\frac{\partial}{\partial Y} +3U\frac{\partial}{\partial U},
\right. \\
{ } & { } & \qquad \left.
 \pm\frac{1}{k+1}\left(Y\pm X^{k+1}\right)\frac{\partial}{\partial X}
+X^k U\frac{\partial}{\partial U}, \;
\pm\frac{U}{k+1}\frac{\partial}{\partial X} +X^k Y\left(Y\pm
X^{k+1}\right)\frac{\partial}{\partial U}, \right.
\\
{ } & { } & \qquad \left. \textcolor{black}{+}\frac{X U}{k+1}\frac{\partial}{\partial
X} +\left(U^2-Y^2\left(Y\pm X^{k+1}\right)\right)
\frac{\partial}{\partial U} \right\rangle_{C_0}.
\end{eqnarray*}
}
\medskip
\par
\textcolor{black}{
Set ${\bf v}_1=2U\frac{\partial}{\partial Y} +
\left({3Y^2\pm 4X^{k+1}Y+X^{2k+2}} \right) \frac{\partial}{\partial U}$, 
${\bf v}_2=\frac{2X}{k+1}\frac{\partial}{\partial X}
+2Y\frac{\partial}{\partial Y} +3U\frac{\partial}{\partial U}$, 
${\bf v}_3=\pm\frac{1}{k+1}\left(Y\pm X^{k+1}\right)\frac{\partial}{\partial X}
+X^k U\frac{\partial}{\partial U}$, 
${\bf v}_4=\pm\frac{U}{k+1}\frac{\partial}{\partial X} +X^k Y\left(Y\pm
X^{k+1}\right)\frac{\partial}{\partial U}$ and 
${\bf v}_5=\frac{X U}{k+1}\frac{\partial}{\partial
X} +\left(U^2-Y^2\left(Y\pm X^{k+1}\right)\right)
\frac{\partial}{\partial U}$.      
Then, we have the following relation:  
\[
-Y{\bf v}_1+U{\bf v}_2\pm X{\bf v}_4={\bf v}_5.    
\]
And, it is easily seen that none of ${\bf v}_1, {\bf v}_2, {\bf v}_3, {\bf v}_4$ 
can be generated by others.      }
\par 
\textcolor{black}{
It is also easily seen that the minimal number of
generators for $Lift(S_k^{\pm})$ is less than or equal to the
minimal number of generators for $Lift(F)$. Thus, the minimal number
of generators for $Lift(S_k^{\pm})$ is less than or equal to $5$. 
It is interesting to observe that the minimal number of generators for $Lift(S_k^{\pm})$ is always $4$.      
It is also interesting to observe that the germs $B_k^{\pm}$, $C_k^{\pm}$
and $F_4$ in Mond's classification (\cite{mond}) also have less than
or equal to 5 generators in the set of liftable vector fields, since
they all admit one-parameter stable unfoldings $\mathcal
A$-equivalent to $F$.
}

\subsection{$Lift(f)$ for $f(x,y)=\{(x,y^3+xy),(x,y^2)\}$}
\label{subsection 6.9} \qquad \\
Let $f=\{f_1,f_2\}$ be the plane to plane bigerm defined by
$f_1(x,y)=(x,y^3+xy)$ and $f_2(x,y)=(x,y^2)$. Consider the
one-parameter stable unfolding $F=\{F_1,F_2\}$ defined by
$$
\begin{cases}
(x,y^3+xy,z)\\
(x,y^2+z,z)
\end{cases}.
$$

It is not hard to see that $Lift(F_1)=\langle
2X\frac{\partial}{\partial X}+3Y\frac{\partial}{\partial
Y},
9Y\frac{\partial}{\partial X}\textcolor{black}{-}2X^2\frac{\partial}{\partial
Y},
\frac{\partial}{\partial Z}\rangle_{C_0}$ 
and $Lift(F_2)=\langle
\frac{\partial}{\partial X}, Y\frac{\partial}{\partial
Y}+Z\frac{\partial}{\partial Z},\frac{\partial}{\partial
Y}+\frac{\partial}{\partial Z}\rangle_{C_0}$ where $(X,Y,Z)$ are the
variables in the target. So $Lift(F)=Lift(F_1)\cap Lift(F_2)=$
$$
\langle (Z-Y)\frac{\partial}{\partial
Z},
2X\frac{\partial}{\partial X}+3Y\frac{\partial}{\partial
Y}+3Z\frac{\partial}{\partial Z},
9Y\frac{\partial}{\partial
X}\textcolor{black}{-}2X^2\frac{\partial}{\partial Y}
\textcolor{black}{-}2X^2\frac{\partial}{\partial
Z}\rangle_{C_0}.$$

To apply Theorem \ref{theorem 2}, we consider $g(x,y,z)=(x,y,z^2)$.
Then $Lift(g)=\langle \frac{\partial}{\partial X},
\frac{\partial}{\partial Y}, Z\frac{\partial}{\partial
Z}\rangle_{C_0}$. 
So $Lift(F)\cap Lift(g)=$
$$
\langle 2X\frac{\partial}{\partial X}+3Y\frac{\partial}{\partial
Y}+3Z\frac{\partial}{\partial Z},
9Y^2\frac{\partial}{\partial X}\textcolor{black}{-}2X^2Y\frac{\partial}{\partial Y}
\textcolor{black}{-}2X^2\textcolor{black}{Z}\frac{\partial}{\partial
Z},$$ 
$$(27YZ\textcolor{black}{+}4X^3)\frac{\partial}{\partial
X}+(\textcolor{black}{-}6X^2Z\textcolor{black}{+}6X^2Y)\frac{\partial}{\partial
Y},(Z^2-YZ)\frac{\partial}{\partial Z}\rangle_{C_0}.$$

And finally 
$Lift(f)=\langle 2X\frac{\partial}{\partial
X}+3Y\frac{\partial}{\partial Y},
9Y^2\frac{\partial}{\partial
X}\textcolor{black}{-}2X^2Y\frac{\partial}{\partial Y} \rangle_{C_0}$.

\medskip
\textcolor{black}{
\begin{remark}\label{remopsu}
Few complete classifications of simple multigerms are known, namely
\cite{kolgushkinsadykov}, \cite{robertathesis}, \cite{hobbskirk} and
\cite{raulcidinharoberta2}. Based on these classifications, it seems
that most simple germs, except for a few cases which can be excluded
by the multiplicity (\cite{nishimura}), admit one-parameter stable
unfoldings. This suggests that the method followed above can be
applied to most simple germs.
\end{remark}}

\section{The case $n>p$}\label{section 7}

Let $\bar{f}:(\mathbb K^n,0)\rightarrow (\mathbb K^p,0)$ be a corank
1 simple germ with $n>p$. By \cite{riegerruas}, $\bar{f}$ is
$\mathcal A$-equivalent to
$$f(x_1,\ldots,x_p,\ldots,x_n)=(x_1,\ldots,x_{p-1},g(x_1,\ldots,x_p)+\sum_{j=p+1}^na_jx_j^2),$$
where $a_j=\pm 1$, $(p+1\leq j\leq n)$. Let $f_0:(\mathbb
K^p,0)\rightarrow (\mathbb K^p,0)$ be the germ such that
$f_0(x_1,\ldots,x_p)=(x_1,\ldots,x_{p-1},g(x_1,\ldots,x_p))$.

\begin{proposition}\label{nbigp}
$Lift(f)=Lift(f_0)$.
\end{proposition}


\par
\bigskip \noindent \underline{\it Proof of Proposition
\ref{nbigp}}\qquad
\par
First suppose $\eta\in Lift(f)$, by definition there exists
$\xi\in\theta_0(n)$ such that $\eta\circ
f(x_1,\ldots,x_n)=df\circ\xi(x_1,\ldots,x_n)$. In particular,
$\eta\circ
f(x_1,\ldots,x_p,0,\ldots,0)=df\circ\xi(x_1,\ldots,x_p,0,\ldots,0)$
and therefore $\eta\circ
f_0(x_1,\ldots,x_p)=df_0\circ\bar{\xi}(x_1,\ldots,x_p)$ where

\[
\bar{\xi}(x_1,\ldots,x_p)= \left( \begin{array}{c}
\xi_1(x_1,\ldots,x_p,0,\ldots,0)\\
\vdots\\
\xi_p(x_1,\ldots,x_p,0,\ldots,0)
\end{array} \right).
\]
Therefore, $\eta\in Lift(f_0)$ and so $Lift(f)\subset Lift(f_0)$.

Now suppose $\eta_0\in Lift(f_0)$, by definition there exists
$\xi_0\in\theta_0(p)$ such that $\eta_0\circ f_0=df_0\circ\xi_0$.
Note that $df_0\circ\xi_0=df\circ\xi$ where
\[
\xi(x_1,\ldots,x_n)= \left( \begin{array}{c}
\xi_{0_1}(x_1,\ldots,x_p)\\
\vdots\\
\xi_{0_p}(x_1,\ldots,x_p)\\
0\\
\vdots\\
0
\end{array} \right).
\]

We have that $\eta_0\circ f(x_1,\ldots,x_n)=$

\[
\left( \begin{array}{c} \eta_{0_1}\circ
f(x_1,\ldots,x_p,0,\ldots,0)\\
\vdots\\
\eta_{0_p}\circ f(x_1,\ldots,x_p,0,\ldots,0)
\end{array} \right) + \left(
\begin{array}{c} \eta_{0_1}\circ
f(x_1,\ldots,x_n)-\eta_{0_1}\circ
f(x_1,\ldots,x_p,0,\ldots,0)\\
\vdots\\
\eta_{0_p}\circ f(x_1,\ldots,x_n)-\eta_{0_p}\circ
f(x_1,\ldots,x_p,0,\ldots,0) \end{array} \right)
\]

The first matrix is equal to $\eta_0\circ
f_0=df_0\circ\xi_0=df\circ\xi$ and, by Hadamard's Lemma
(\cite{gibson}), there exist functions
$\bar{\xi}_{p+1},\ldots,\bar{\xi}_n$ such that the second matrix is
equal to
\[
df\circ \left( \begin{array}{c} \eta_{0_1}\circ
f(x_1,\ldots,x_n)-\eta_{0_1}\circ
f(x_1,\ldots,x_p,0,\ldots,0)\\
\vdots\\
\eta_{0_p}\circ f(x_1,\ldots,x_n)-\eta_{0_p}\circ
f(x_1,\ldots,x_p,0,\ldots,0)\\
0\\
\bar{\xi}_{p+1}\\
\vdots\\
\bar{\xi}_n
\end{array} \right).
\]

So $\eta_0\in Lift(f)$ and the proposition is proved.
  \hfill Q.E.D.

The proposition holds for multigerms too since the above proof can
be repeated for each branch. Thus, we can obtain the following

\begin{example}
Let $f=\{f_1,f_2\}:(\mathbb K^n,\{0,0\})\rightarrow (\mathbb
K^2,0)$, $n>2$, be the bigerm defined by
$$
\begin{cases}
f_1(x,y,u_1,\ldots,u_{n-2})=(x,y^3+xy+\sum_{i=1}^{n-2}a_iu_i^2)\\
f_2(x,y,u_1,\ldots,u_{n-2})=(x,y^2+\sum_{i=1}^{n-2}b_iu_i^2)
\end{cases}
$$
where $a_i=\pm 1$ and $b_i=\pm 1$, $(1\leq i\leq n-2)$. Then
$Lift(f)=\langle 2X\frac{\partial}{\partial
X}+3Y\frac{\partial}{\partial Y},9Y^2\frac{\partial}{\partial
X}+2X^2Y\frac{\partial}{\partial Y} \rangle_{C_0}$.
\end{example}
%
%
\textcolor{black}
{
\section*{Acknowledgements}
The authors would like to thank Mohammed Salim Jbara Al-Bahadeli for his pointing out careless mistakes in calculations for concrete liftable vector fields       
}
%
%
%

\end{document}